\documentclass{amsart} 
\usepackage{amssymb}
\usepackage{amsmath}
\usepackage{amsfonts}
\usepackage{amsthm}
\usepackage{color}
\begin{document} 
\input epsf.sty

\newtheorem{lem}{Lemma}[section]
\newtheorem{prop}[lem]{Proposition}
\newtheorem{cor}[lem]{Corollary}
\newtheorem{thm}[lem]{Theorem}

\date{}


\address{F. Pakovich, Department of Mathematics,
Ben Gurion University,
P.O.B. 653,
Beer Sheva 84105, Israel}
\email{pakovich@math.bgu.ac.il}
\address{M. Muzychuk, Department of Mathematics, Netanya Academic College, Kibbutz Galuyot St. 16, Netanya 42365, Israel}
\email{muzy@netanya.ac.il}

\thanks{The research of the first author was supported by the ISF grant no. 979/07}

\begin{abstract} In this paper we give a solution of the following ``polynomial moment problem'' which arose about ten years ago in connection with Poincar\'e's center-focus problem: for a given polynomial $P(z)$ to describe polynomials $q(z)$ orthogonal to all powers of $P(z)$ on a segment $[a,b]$. 
\end{abstract}

\title{Solution of the polynomial moment problem
}  

\author{F. Pakovich, M. Muzychuk}

\maketitle

\date{}

\newcommand{\new}[1]{{\emph {#1}}}
\newcommand{\aut}{{\sf Aut}}
\newcommand{\sym}{{\sf Sym}}
\newcommand{\orb}{{\sf Orb}_2}
\newcommand{\A}{\mathbb{A}}
\newcommand{\C}{\mathbb{C}}
\newcommand{\F}{\mathbb{F}}
\newcommand{\N}{\mathbb{N}}
\newcommand{\Q}{\mathbb{Q}}
\newcommand{\Z}{\mathbb{Z}}

\newcommand{\cA}{\mathcal{A}}
\newcommand{\cB}{\mathcal{B}}
\newcommand{\cR}{\mathcal{R}}
\newcommand{\cV}{\mathcal{V}}
\newcommand{\cS}{\mathcal{S}}
\newcommand{\cX}{\mathcal{X}}
\newcommand{\wh}[1]{\widehat{#1}}
\newcommand{\sig}{\sigma}
\newcommand{\lm}{\lambda}
\newcommand{\im}{{\sf Im}}
\newcommand{\Ker}{\mathrm{Ker}}
\newcommand{\wt}[1]{\widetilde {#1}}
\renewcommand{\a}{\tilde{a}}
\newcommand{\tR}{\tilde{R}}
\newcommand{\irr}{{\sf Irr}}
\newcommand{\tr}{\mathrm{Tr}}
\newcommand{\un}{\underline}
\newcommand{\fX}{\mathfrak X}
\newcommand{\Trace}[1]{\mathaccent23{#1}}
\newcommand{\New}[1]{{\sl #1}}
\newcommand{\und}[1]{\underline{#1}}
\newcommand{\sg}[1]{\langle #1\rangle}
\newcommand{\Bsets}[1]{{\sf Basic}(#1)}
\newcommand{\eps}{\epsilon}
\newcommand{\misha}[1]{\textcolor{red}{#1}}

\def\be{\begin{equation}}
\def\ee{\end{equation}}
\def\bs{$\square$ \vskip 0.2cm}
\def\d{{\rm d}} 
\def\D{{\rm D}} 
\def\I{{\rm I}} 
\def\C{{\mathbb C}} 
\def\N{{\mathbb N}} 
\def\P{{\mathbb P}}
\def\Q{{\mathbb Q}}
\def\Z{{\mathbb Z}}
\def\R{{\mathbb R}} 
\def\ord{{\rm ord}}
\def\qed{$\ \ \Box$ \vskip 0.2cm}

\def\e{\eqref}
\def\phi{{\varphi}}
\def\v{{\varepsilon}} 
\def\deg{{\rm deg\,}} 
\def\Det{{\rm Det}}
\def\dim{{\rm dim\,}} 
\def\Ker{{\rm Ker\,}} 
\def\Gal{{\rm Gal\,}}
\def\St{{\rm St\,}} 
\def\exp{{\rm exp\,}} 
\def\cos{{\rm cos\,}}
\def\card{{\rm card\,}} 
\def\diag{{\rm diag\,}} 
\def\GCD{{\rm GCD }}
\def\LCM{{\rm LCM }}
\def\mod{{\rm mod\ }}

\def\bp{\begin{proposition}}
\def\ep{\end{proposition}}
\def\bt{\begin{theorem}}
\def\et{\end{theorem}}
\def\be{\begin{equation}}
\def\l{\label}
\def\ee{\end{equation}}
\def\bl{\begin{lemma}}
\def\el{\end{lemma}}
\def\bc{\begin{corollary}}
\def\ec{\end{corollary}}
\def\pr{\noindent{\it Proof. }}
\def\note{\noindent{\bf Note. }}
\def\bd{\begin{definition}}
\def\ed{\end{definition}}

\newtheorem{theorem}{Theorem}[section]
\newtheorem{lemma}[theorem]{Lemma}
\newtheorem{definition}[theorem]{Definition}
\newtheorem{corollary}[theorem]{Corollary}
\newtheorem{proposition}[theorem]{Proposition}

\section{Introduction}

In this paper we solve the following
``polynomial moment problem'':
{\it for given $P(z)\in \C[z]$ and distinct $a,b\in \C$ 
to describe $q(z)\in \C[z]$ such that
\be \l{1}
\int^b_a P^i(z)q(z)\d z=0 
\ee
for all $i\geq 0$.}

The polynomial moment problem was posed in the 
series of papers \cite{bfy1}-\cite{bfy4}
in connection with the
center problem for the Abel differential equation 
\be \l{ab}
\frac{\d y}{\d z}=p(z)y^2+q(z)y^3.
\ee
with polynomial
coefficients $p(z), q(z)$ in the complex domain.
For given $a,b\in \C$ the center problem for the Abel equation
is to find necessary and sufficient conditions on $p(z),q(z)$ which imply
the equality $y(b)=y(a)$ for any solution $y(z)$ of \eqref{ab} with $y(a)$ small enough.
This problem is closely related to the classical 
Center-Focus problem of Poincar\'e and
has been studied in many recent papers 
(see e.g. \cite{bby}-\cite{c}, \cite{y1}).

The center problem for the Abel equation is connected with the 
polynomial moment problem in several ways.
For example, it was shown in \cite{bfy3}
that for the parametric version  
$$
\frac{\d y}{\d z}=p(z)y^2+\varepsilon q(z)y^3
$$
of \eqref{ab} the ``infinitesimal'' center conditions with respect to $\varepsilon$ 
reduce to moment equations \eqref{1} with $P(z)=\int p(z) \d z.$
On the other hand, it was shown in \cite{bry} that ``at infinity'' (under an appropriate projectivization of the 
parameter space) the system of equations on the
coefficients of $q(z)$, describing the center set of \eqref{ab} for fixed $p(z)$,
also reduces to \eqref{1}.
Many other results concerning connections between 
the center problem and the polynomial moment problem can be found in \cite{bry}.
These results convince that a thorough 
description of solutions of system \eqref{1} is an important step 
in the understanding of the center problem for the Abel equation.

There exists a natural condition on $P(z)$ and $Q(z)=\int q(z)\d z$
which reduces equations \eqref{1}, \eqref{ab} to similar equations
with respect to polynomials of smaller degrees. Namely, suppose that
there exist polynomials $\tilde P(z),$ $\tilde Q(z),$ $W(z)$ with $\deg W(z)>1$ such that
\be \l{2}
P(z)=\tilde P(W(z)), \ \ \ \ \ 
\ Q(z)=\tilde Q(W(z)).
\ee
Then after the change of variable $w=W(z)$
equations \eqref{1} transform to the equations  
\be \l{zam}
\int^{W(b)}_{W(a)} \tilde P^i(w)\tilde Q^{\prime}(w)\d w=0, 
\ee
while equation \eqref{ab} transforms to the equation 
\be \l{ab1}
\frac{\d\tilde y}{\d w}=\tilde P^{\prime} (w)\tilde y^2+\tilde Q^{\prime}(w)\tilde y^3.
\ee

Furthermore, if the polynomial $W(z)$ in \eqref{2} 
satisfies the equality 
\be \l{w} W(a)=W(b), 
\ee 
then the Cauchy theorem implies that the polynomial $\tilde Q^{\prime}(w)$ is a solution of system \eqref{zam} and hence the polynomial $q(z)=Q^{\prime}(z)$ is a solution of system \eqref{1}.
Similarly, since any solution $y(z)$ of equation \eqref{ab} is the pull-back 
\be \l{pb} y(z)=\tilde y(W(z))
\ee
of a solution $\tilde y(w)$ 
of equation \eqref{ab1}, if $W(z)$ satisfies \eqref{w} then equation
\eqref{ab} has a center. 
This justifies the following definition:
a center for equation
\eqref{ab} or a solution of system \eqref{1} is called
{\it reducible} if there exist polynomials 
$\tilde
P(z),$ $\tilde Q(z),$ $W(z)$ such that
conditions 
\eqref{2}, \eqref{w} hold.
The main conjecture concerning the center problem for the Abel equation
(``the composition conjecture for the Abel equation"), supported by the results obtained 
in the papers cited above, states 
that any center for the Abel equation is reducible 
(see \cite{bry}
and the bibliography there).

By analogy with the composition conjecture it was suggested 
(``the composition conjecture for the polynomial moment problem") 
that, under the additional assumption $P(a)=P(b),$
any solution of (1) is reducible.
This conjecture was shown to be 
true in many cases. For instance, 
if $a,b$ are not critical points of $P(z)$   
(\cite{c}), if $P(z)$ is indecomposable (\cite{pa2}), and in some other 
special cases (see e. g. \cite{bfy3}, \cite{pp}, \cite{pry}, \cite{ro}).
Nevertheless, 
in general the composition conjecture for the polynomial moment problem 
fails to be true. 
Namely, it was shown in \cite{pa1} that if $P(z)$ has several ``compositional right factors'' $W(z)$ such that $W(a)=W(b)$, then
it may happen that the sum of reducible solutions corresponding to these factors is a non-reducible solution.
 
It was conjectured in \cite{pa4} that actually {\it any} non-reducible solution of \eqref{1} is a sum of reducible ones.
Since compositional factors $W(z)$ of a polynomial $P(z)$ can be defined explicitly,
such a description of non-reducible solutions of \eqref{1}
would be very helpful, especially   
for 
applications to the 
Abel equation
(cf. \cite{bry}).
However, until now
this conjecture was verified only in a single special case 
(see \cite{pa3}). 

Meanwhile, another necessary and sufficient condition 
for a polynomial $q(z)$ to be a solution of \eqref{1} was constructed in \cite{pp}. 
Namely, it was shown in \cite{pp} that there exists a finite system of equations 
\be \l{su}
\sum_{i=1}^nf_{s,i}Q(P^{-1}_{i}(z))=0, \ \ \ \ \ \ f_{s,i} \in \Z,  \ \ \ \ \ \ 1\leq s\leq k, 
\ee 
where $Q(z)=\int q(z)\d z$ and 
$P^{-1}_{i}(z),$ 
$1\leq i \leq n,$ are branches of the algebraic function  
$P^{-1}(z)$, 
such that \eqref{1} holds if and only if \eqref{su} holds.
Moreover, this system was constructed 
explicitly with the use of a special planar tree $\lambda_P$ 
which represents the monodromy group $G_P$ of the algebraic function $P^{-1}(z)$ in a combinatorial way.
By construction, points $a,b$ are vertices of $\lambda_P$ and system 
\eqref{su} reflects the combinatorics of the path connecting $a$ and $b$ on $\lambda_P$.

A finite system of equations 
\eqref{su} 
is more convenient for a study than initial infinite system of equations \eqref{1}.
In particular, in many cases the analysis of \eqref{su} permits to conclude that for given $P(z),a,b$ any solution 
of \eqref{1} is reducible (see \cite{pp}).
In this paper we develop
necessary algebraic and analytic techniques which allow us to describe 
solutions of \eqref{su} in the general case
and to prove that any solution of \eqref{1} is a sum of reducible ones.
So, our main result is the following theorem.

\bt \l{os} A non-zero polynomial $q(z)$ is a solution of system (1) 
if and only if $Q(z)=\int q(z)\d z$ can be represented as a sum of
polynomials $Q_j(z)$ such that
\be \l{cc}
P(z)=\tilde P_j(W_j(z)), \ \ \
Q_j(z)=\tilde Q_j(W_j(z)), \ \ \ {\it and} \ \ \ W_j(a)=W_j(b)
\ee
for some polynomials $\tilde P_j(z), \tilde Q_j(z), W_j(z)$.
\et

Note that since conditions of the theorem 
impose no restrictions on
the values of $P(z)$ at the points $a,b$ the theorem implies in particular 
that non-zero solutions of \eqref{1} exist if 
and only if the equality
$P(a)=P(b)$ holds. Indeed, if $P(a)=P(b)$ then
for any $\tilde Q(z)\in \C[z]$ the polynomial $Q(z)=\tilde Q(P(z))$ is a solution of \eqref{1} since we can set
$W(z)=P(z)$ in \eqref{2}, \eqref{w}. On the other hand, if $Q(z)$ is a solution of \eqref{1} then equalities \eqref{cc} imply that $P(a)=P(b)$. 

The paper is organized as follows. 
In the second section we give a detailed 
account of definitions and previous results related to the polynomial moment 
problem.   
In particular, starting from system \eqref{su}, we introduce a linear subspace $M_{P,a,b}$ of 
$\mathbb Q^n$ 
generated by the vectors 
$$(f_{s,\sigma(1)}, f_{s,\sigma(2)},\, ...\,, f_{s,\sigma(n)}), \ \ \ \sigma\in G_P, \ \ \ 1\leq s \leq k,$$
and study its basic properties. 

It follows from the definition that $M_{P,a,b}$ is invariant 
with respect to the permutation matrix representation of 
the group 
$G_P$.
In the third section of the paper, written entirely in the framework of the group theory,
we describe a general structure of such subspaces.
More generally, 
we describe subspaces of $\mathbb Q^n$ 
invariant with respect to the 
permutation matrix representation of a permutation group $G$ 
of degree $n$, containing a cycle of length $n$. Roughly speaking, we show that the structure of invariant subspaces of
$\mathbb Q^n$ for such $G$ depends on imprimitivity systems of $G$ only. 
We believe that this result is new and interesting by itself. 

Finally, in the fourth section, using the description of $G_P$-invariant subspaces of $\mathbb Q^n$  
and results and techniques of \cite{pp}, we prove Theorem \ref{os}.

\vskip 0.2cm
\noindent{\bf Acknowledgments}. The authors would
like to thank C. Christopher, J. P. Fran\-\c{c}oise, L. Gavrilov, G. Jones, M. Klin, Y. Yomdin, and W. Zhao 
for discussions of different questions
related to the subject of this paper. The authors also are grateful to the anonymous 
referee for valuable comments and suggestions.

\section{Preliminaries
} In this section we collect basic definitions and results concerning the polynomial moment problem.
In order to make the paper self-contained we outline
proofs of main statements.

\subsection{Criterion for $\hat H(t)\equiv 0$}
For $P(z),Q(z)\in \C[z]$ and 
a path $\Gamma_{a,b}\subset \C$ connecting different points $a,b$ of $\C$  
let
$H(t)=H(P,Q,\Gamma_{a,b},t)$ be a function defined on $\C\P^1\setminus P(\Gamma_{a,b})$ by the integral
\be \l{cau}
H(t)= \int_{\Gamma_{a,b}} \frac{Q(z)P^{\prime}(z)dz}{P(z)-t}.
\ee
Notice that although integral \eqref{cau} depends on $\Gamma_{a,b}$ the Cauchy theorem implies that if $\tilde \Gamma_{a,b}\subset \C$ is another path connecting $a$ and $b$,
then for all $t$ close enough to infinity the equality $$H(P,Q,\tilde \Gamma_{a,b},t)=H(P,Q,\Gamma_{a,b},t)$$
holds. Therefore, the Taylor expansion of $H(t)$ at infinity and the corresponding germ $\hat H(t)$ do
not depend on the choice of $\Gamma_{a,b}.$

After the change of variable
$z\rightarrow P(z)$
integral \eqref{cau} transforms to the Cauchy type integral
\be \l{ci} H(t)=
\int_{\gamma}\frac{g(z)dz}
{z-t}\,,
\ee where $\gamma=P(\Gamma_{a,b})$ and
$g(z)$ is obtained by the analytic continuation
along $\gamma$ of a germ of the algebraic function $Q(P^{-1}(z)).$ 
Clearly, integral representation \eqref{ci}
defines an analytic function in each domain of the complement
of $\gamma$ in $\C\P^1$. Notice that for any choice of $\Gamma_{a,b}$
the function defined in the domain containing infinity is the analytic continuation of the germ $\hat H(t).$

\bl [\cite{pp}]\l{lem} Assume that $P(z),$ $q(z)\in \C[z]$ 
and $a,b\in \C,$ $a\neq b,$ satisfy 
\be \l{ux} \int_{\Gamma_{a,b}} P^i(z)q(z)\d z=0, \ \ \ i\geq 0, \ee
and let $Q(z)$ be a polynomial defined by the equalities \be \l{urur} Q(z)=\int q(z) \d z, \ \ \ Q(a)=0.\ee
Then for the germ $\hat H(t)$ defined near infinity by integral \eqref{cau}  
the equality $\hat H(t)\equiv 0$ holds.
\el

\pr Indeed, for all $i\geq 1$ by integration by parts we have:
\be \l{ururu} 
\int_{\Gamma_{a,b}} P^i(z)q(z)\d z=P^i(b)Q(b)-P^i(a)Q(a)
-i\int_{\Gamma_{a,b}} P^{i-1}(z)Q(z)P^{\prime}(z)\d z.\ee
Furthermore, $Q(a)=0$ implies $Q(b)=0$ in view of \eqref{ux} taken for $i=0$. Therefore, if \eqref{ux} holds then  
all the integrals appearing in the right part of \eqref{ururu} vanish. 
On the other hand, these integrals 
are coefficients of the Taylor expansion of $-\hat H(t)$ at infinity.
\qed

Lemma \ref{lem} shows that the polynomial moment problem reduces to the problem of 
finding conditions on $Q(z)$ under which the equality $\hat H(t)\equiv 0$ holds.
On the other hand, we will show below (Corollary \ref{zzaa}) that if $\hat H(t)\equiv 0$ holds for some polynomial $Q(z)$ then 
\eqref{ux} holds for $q(z)=Q'(z).$ 
A condition of a general character for $\hat H(t)$ to vanish was given in the paper \cite{pry}
in the context of the theory of Cauchy type integrals of algebraic functions.
Subsequently, in the paper \cite{pp} was proposed a construction
which permits to obtain conditions for
the vanishing of $\hat H(t)$ in a very explicit form.
Briefly, the idea of \cite{pp} is to choose the integration path $\Gamma_{a,b}$
in such a way that its image under the mapping $P(z)\,:\, \C\P^1\rightarrow \C\P^1$ does not divide the Riemann sphere.

The construction of the paper \cite{pp} uses
a special graph $\lambda_P,$
embedded into the Riemann sphere, defined as follows
(see \cite{pp}).
Let $S$ be
a ``star'' joining a non-critical value $c$ of a polynomial $P(z)$ of degree $n$ with all its {\it finite} critical values $c_1,c_2, ... ,c_k$ by non intersecting oriented arcs
$\gamma_1, \gamma_2, ... ,\gamma_k$. Define $\lambda_P$
as a preimage of $S$ under the map $P(z)\,:\, \C\P^1\rightarrow \C\P^1$ (see Fig. 1).
\begin{figure}[ht]
\epsfxsize=12truecm
\centerline{\epsffile{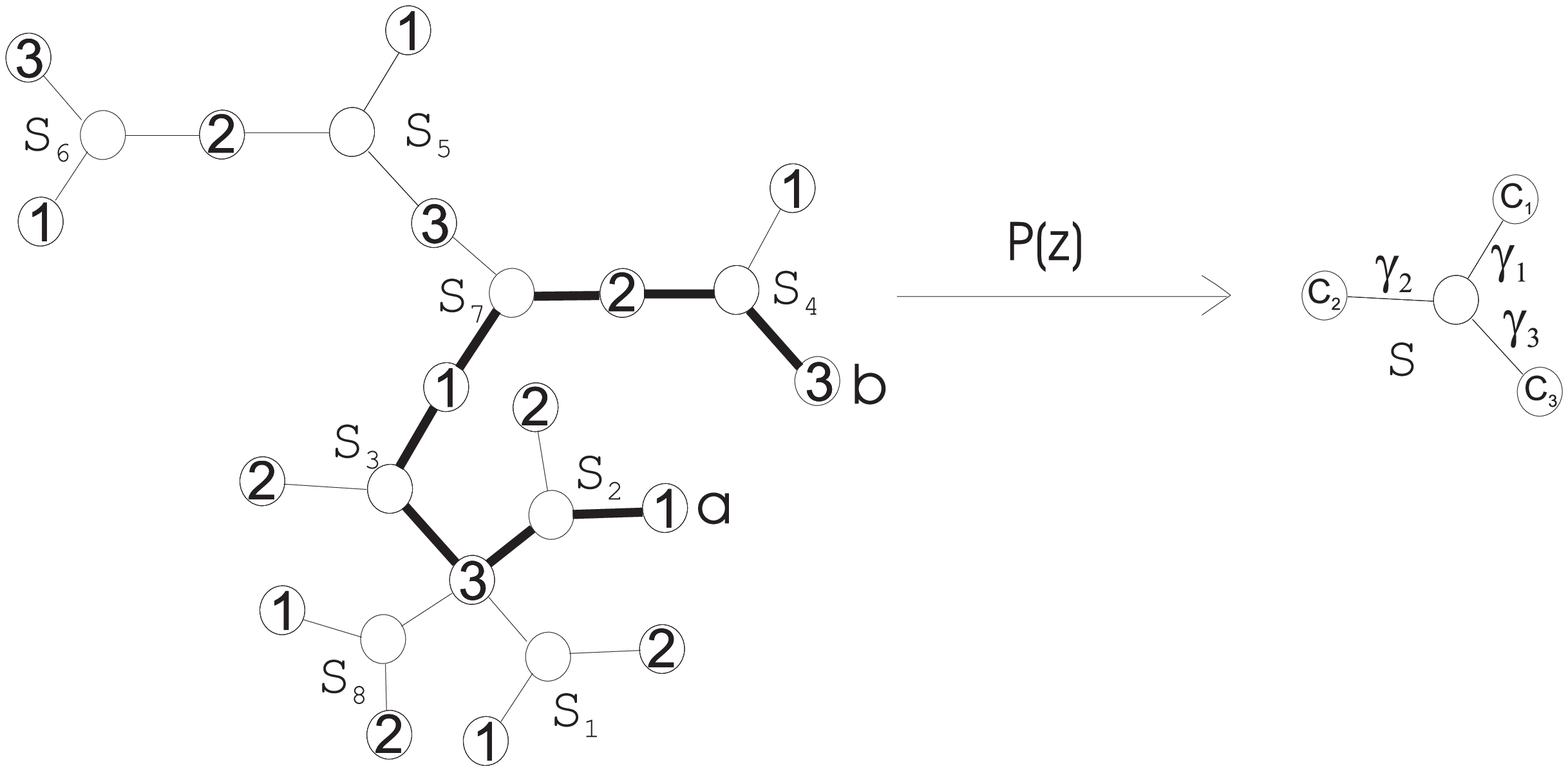}}
\smallskip
\centerline{Figure 1}
\end{figure}
More precisely, define
vertices of $\lambda_P$ as
preimages
of the points $c$ and $c_s,$ $1\leq s \leq k,$
and edges of $\lambda_P$ as preimages of the arcs
$\gamma_s,$ $1 \leq s \leq k.$
Furthermore, for each $s,$ $1\leq s \leq k,$ mark vertices of $\lambda_P$ which are preimages of the point $c_s$ by the number $s.$ Finally, define a {\it star} of $\lambda_P$ as a subset of edges of $\lambda_P$ consisting of edges adjacent to some
non-marked vertex.

By construction, the restriction of $P(z)$ on $\C\P^1\setminus \lambda_P$ is a covering of the topological punctured disk $\C\P^1\setminus \{S\cup \infty\}$ and therefore $\C\P^1\setminus \lambda_P$ is a disjointed union of punctured disks (see e.g. \cite{fors}). Moreover, since the preimage of infinity under $P(z)$ consists of a unique point, $\C\P^1\setminus \lambda_P$ consists of a unique disk and hence the graph $\lambda_P$ is a tree.

Set $C=\{c_1,c_2,
... ,c_k\}$ and let $U\subset \C$ be a simply connected domain
such that $S\setminus C\subset U$ and $U\cap C=\emptyset$. Then in $U$
there exist $n$ single-valued analytic 
branches of the algebraic function $P^{-1}(z)$ inverse to $P(z).$ We will denote these branches by $P^{-1}_i(z),$ $1\leq i \leq n.$
The stars
of $\lambda_P$ may be naturally identified with branches of $P^{-1}(z)$ in $U$ as follows:
to the branch $P^{-1}_i(z),$ $1\leq i \leq n,$
corresponds the star $S_i,$ $1\leq i \leq n,$
such that $P^{-1}_i(z)$ maps bijectively the interior of $S$ to the interior of $S_i.$

Under the analytic continuation along a closed curve  
the set $P^{-1}_i(z),$ $1\leq i \leq n,$ transforms 
to itself. This induces a homomorphism \be \l{homo} \pi_1(\C\P^1\setminus \{C\cup \infty\},c)\rightarrow S_n\ee whose 
image $G_P$ is called the monodromy
group of $P(z)$. 
Notice that if $\omega_{\infty}$ and $\omega_i$, $1\leq i \leq k,$ are loops around $\infty$ and
$c_i$, $1\leq i \leq k,$ respectively, such that $\omega_1\omega_2\dots \omega_k\omega_{\infty}=1$ in 
$\pi_1(\C\P^1\setminus  \{C\cup \infty\},c),$ then the elements $g_i$, $1\leq i \leq k,$ of $G_P$,
which are the images of $\omega_i$, $1\leq i \leq k,$ under homomorphism \eqref{homo},
generate $G_P$ and satisfy the equality $g_1g_2\dots g_kg_{\infty}=1,$ where $g_{\infty}$ is the element of $G_P$ which is the image of $\omega_{\infty}$.

Having in mind the identification of
the set of stars of $\lambda_P$ with the set of branches of $P^{-1}(z),$
the permutation $g_s,$ $1\leq s \leq k,$ can be
identified with a permutation $\hat g_s,$
$1\leq s \leq k,$
acting
on the set of stars of $\lambda_P$ in the following way: $\hat g_s$ sends the
star $S_i,$ $1\leq i \leq n,$ to the ``next'' star under a counterclockwise
rotation around the vertex of $S_i$ colored by the $s$th color. For example,
for the tree shown on Fig. 1 we have: $g_1=(1)(2)(37)(4)(5)(6)(8),$
$g_2=(1)(2)(3)(47)(56)(8),$ $g_3=(1238)(4)(57)(6).$ Notice that since $P(z)$ is a polynomial, the permutation
$g_{\infty}$ is a cycle of
length $n.$ We always will assume that the numeration of branches of $P^{-1}(z)$ in $U$ is chosen in such a way
that $g_{\infty}$ coincides with the cycle $(1\,2\,...\,n)$. Clearly, such a numeration is
defined uniquely up to a choice of $P_1^{-1}(z)$.

The tree constructed above is known under the name of ``constellation'' or ``cactus'' and is closely related to
what is called a ``dessin d'enfant" (see \cite{lz}
for further details and other versions of this construction).
Notice that the Riemann existence theorem implies that
a polynomial $P(z)$ is defined by $c_1,c_2, ... ,c_k$
and $\lambda_P$ up to a composition $P(z)\rightarrow P(\mu (z)),$ where $\mu(z)$ is a linear function.

It follows from the definition that the points $a$ and $b$ are vertices of $\lambda_P$ if and only if $P(a)$ and $P(b)$
are critical values of $P(z)$. For our purposes however it is more convenient
to define the tree $\lambda_P$ so that
the points $a,b$ always would be its vertices.
So, in the case when $P(a)$ or $P(b)$ (or both of them) is not a
critical value of $P(z)$ we modify the construction as follows. Define $c_1,c_2, ... ,c_{k}$ as the set of all finite critical values
of $P(z)$ {\it supplemented} by $P(a)$ or $P(b)$ (or by both of them), and set as above
$\lambda_P=P^{-1}\{S\},$ where $S$ is a star
connecting $c$ with $c_1,c_2, ... ,c_{k}$
(we suppose that $c$ is chosen
distinct from $P(a),P(b)$). Clearly, $\lambda_P$ is still a tree and
the points $a,b$ are vertices of $\lambda_P.$

Since
$\lambda_P$ is connected and has no cycles there exists a unique oriented path $\mu_{a,b}
\subset \lambda_P$
joining the point $a$ with the point $b$.
Furthermore, it follows from the definition of $\lambda_P$ that if we set
$\Gamma_{a,b}=\mu_{a,b}$ then after the change of variable $z\rightarrow P(z)$
integral \eqref{cau}
reduces to the sum of integrals
\be \l{f}
H(t)=\sum_{s=1}^{k}\int_{\gamma_s}\frac{\phi_s(z)}{z-t}\, \d z,
\ee
where each $\phi_s(z),$ $1\leq s \leq k,$ is a linear combination
of the functions $Q(P^{-1}_{i}(z)),$ $1\leq i \leq n,$ in $U$.
Namely, \be \l{piz}
\phi_s(z)=\sum_{i=1}^nf_{s,i}Q(P^{-1}_{i}(z)),
\ee
where $f_{s,i}\neq 0$ if and only if the path
$\mu_{a,b}$ goes through the star
$S_i$ across its $s$-vertex.
Furthermore, if when going along $\mu_{a,b}$
the $s$-vertex of $S_i$
is followed by the center of $S_i$ then $f_{s,i}=-1$ otherwise
$f_{s,i}=1$.
For example, for the graph $\lambda_P$ shown on Fig. 1 and
the path $\mu_{a,b}\subset \lambda_P$ pictured by the fat
line we have:
\begin{gather*}
\phi_1(z)=-Q(P^{-1}_{2}(z))+Q(P^{-1}_{3}(z))- Q(P^{-1}_{7}(z)),\\
\phi_2(z)=Q(P^{-1}_{7}(z))-Q(P^{-1}_{4}(z)),\\
\phi_3(z)=Q(P^{-1}_{2}(z))-Q(P^{-1}_{3}(z))+ Q(P^{-1}_{4}(z)).
\end{gather*}

Notice that the number $k$ in \eqref{f} coincides with the number of critical values $s$ of $P(z)$ such that
the path $\Gamma_{a,b}$ passes through at least one vertex
colored by the $s$-th color.
Note also that
equations \e{piz}
are linearly dependent. Indeed,
for each $i,$ $1\leq i \leq n,$ such that there exists an index $s,$
$1\leq s
\leq k,$ with $f_{s,i}\neq 0$ there exist exactly two such
indices $s_1,s_2$, and $c_{s_1,i}=-c_{s_2,i}.$ Therefore, the equality
$$\sum_{s=1}^{k}\phi_s(t)=0$$
holds in $U$.

\bt [\cite{pp}]\l{t1}
Let $P(z),Q(z)\in \C[z]$ and $a,b\in \C,$ $a\neq b.$
Then $\hat H(t)\equiv 0$ if and only if
$\phi_s(z)\equiv 0$ for any $s,$ $1\leq s \leq k.$
\et

\pr Formula \eqref{f} defines the analytic continuation of
$\hat H(t)$ to the domain
$\C\P^1\setminus S$. In particular, $\hat H(t)\equiv 0$
if and only if $H(t)\equiv 0$ in $\C\P^1\setminus S$.
On the other hand, by the well-known boundary property of
Cauchy type integrals (see e.g. \cite{mus}), for any $s,$ $1\leq s \leq k,$ and any interior point $z_0$ of $\gamma_s$
we have:
\be \l{cuu}2\pi \sqrt{-1} \,\phi_s(z_0)=
\lim_{t \to z_0}\!\!\!\!\,^+H(t)-\lim_{t \to z_0}
\!\!\!\!\,^-H(t),
\ee
where the limits are taken when $t$ approaches $z_0$ from 
the ``right''  (resp. ``left'') side of $\gamma_s$. Therefore, if $H(t)\equiv 0$ in $\C\P^1\setminus S$,
then the limits in \eqref{cuu} equal zero and hence
$\phi_s(z)\equiv 0$ for any $s,$ $1\leq s \leq k.$

Finally,
if \be \l{sss}\phi_s(z)\equiv 0, \ \ \ \ 1\leq s \leq k,\ee then
it follows directly from formula \e{f} that
$H(t)\equiv 0$. \qed

\subsection{Subspace $M_{P,a,b}$}
For any element $\sigma\in G_P$ the
equality $\phi_s(z)=0,$ $1\leq s \leq k,$ implies by the analytic continuation
the equality
$$
\sum_{i=1}^nf_{s,i}Q(P^{-1}_{\sigma(i)}(z))=0.
$$
Therefore, replacing $\sigma$ by $\sigma^{-1}$ we see that Theorem \ref{t1} implies that
$\hat H(t)\equiv 0$ if and only if
for any $\sigma\in G_P$ and $s,$ $1\leq s \leq k,$
the equality
$$
\sum_{i=1}^nf_{s,\sigma(i)}Q(P^{-1}_{i}(z))=0
$$ holds.

Denote by $M_{P,a,b}$ the subspace of $\mathbb Q^n$ generated by the vectors
$$(f_{s,\sigma(1)}, f_{s,\sigma(2)},\, ...\,, f_{s,\sigma(n)}), \ \ \ 1\leq s \leq k, \ \ \ \sigma\in G_P.$$
Abusing the notation we usually will not distinguish an element
of $M_{P,a,b}$ and the corresponding equation connecting branches
of $ Q(P^{-1}(z))$.
For example, instead of using
the notation
\be \l{gi} (0,0,\, ... \,, 1, \, ...\,, 0,0,\, ... \,, -1,\, ... \,,0,0) \ee
for an element of $M_{P,a,b}$ we simply will use the equality
\be \l{compp} Q(P^{-1}_{i_1}(z))=Q(P^{-1}_{i_2}(z)),\ee
for corresponding $i_1\neq i_2,$ $1\leq i_1,i_2\leq n.$

Equality \eqref{compp}
is the simplest example of the equality
$\phi_s(z)=0,$ $1\leq s \leq k,$ 
and is equivalent to the statement that $P(z)$ and $Q(z)$ have a non-trivial ``compositional right factor''
(cf. \cite{c}, \cite{ro}, \cite{pa2},
\cite{pry}, \cite{pp}).

\bl \l{lll} Let $P(z),Q(z)\in \C[z]$. Then
the equalities
\be \l{comp} P(z)=\tilde P(W(z)), \ \ \ \ \ \ Q(z)=\tilde Q(W(z))\ee
hold for some $\tilde P(z), \tilde Q(z), W(z) \in \C[z]$ with $\deg W(z)>1$ if and only if equality \eqref{compp} holds for some $i_1\neq i_2,$ $1\leq i_1,i_2\leq n.$
Furthermore, $Q(z)=\tilde Q(P(z))$ for some $\tilde Q(z)\in \C[z]$ if and only if 
all the functions
$Q(P^{-1}_{i}(z)),$ $1\leq i \leq n,$ are equal between themselves.

\el
\pr Let $d(Q(P^{-1}))$ be a number of {\it different} functions in the collection of functions $Q(P^{-1}_i(z)),$ $1\leq i \leq n.$
Since any algebraic relation over $\C$ between $Q(p^{-1}(z))$ and $z,$ where $p^{-1}(z)$ is a branch of the algebraic function
$P^{-1}(z)$ in $U$, supplies
an algebraic relation between $Q(z)$ and $P(z)$ and vice versa,
we have:
$$d(Q(P^{-1}))=[\C(Q,P):\C(P)]=[\C(z):\C(P)]/[\C(z):\C(Q,P)]=$$ $$=n/[\C(z):\C(Q,P)].$$
Therefore, \be \l{degrr} [\C(z):\C(Q,P)]=n/d(Q(P^{-1})).\ee
It follows now from the L\"{u}roth theorem that $d(Q(P^{-1}))<n$ if and only if   
\eqref{comp} holds for some {\it rational} functions $\tilde P(z),$ $\tilde Q(z),$ $W(z)$ with
$\deg W(z)>1.$ Furthermore, if $d(Q(P^{-1}))=1$ then 
\eqref{degrr} implies that $Q(z)=\tilde Q(P(z))$ for some $\tilde Q(z)\in \C(z).$

Observe now that, since $P(z),$ $Q(z)$ are polynomials, without loss of generality
we may assume that $\C(Q,P)=\C(W)$ for some {\it polynomial} $W(z).$
Indeed, since $P(z)$ is a polynomial the equality $P(z)=U(V(z)),$ where $U(z),$ $V(z)$ are rational functions,
implies that $U(z)$ has a unique pole and that the preimage
of this pole under $V(z)$ consists of infinity only. This implies that $V(z)=\mu(W(z))$
for some polynomial $W(z)$ and M\"obius transformation $\mu(z)$, and it is clear that the fields $\C(V(z))$ and $\C(W(z))$ coincide. Finally, if $W(z)$ is a polynomial then obviously $\tilde P(z),$ $\tilde Q(z)$ also are polynomials. \qed

Since \eqref{comp} implies that
$$
\int^b_a P^i(z)q(z)\d z=
\int^{W(b)}_{W(a)} \tilde P^i(W)\tilde Q^{\prime}(W)\d W,
$$
Lemma \ref{lll} shows that
if the subspace $M_{P,a,b}$ contains an element of the form \eqref{compp},
then any solution $q(z)$ of the polynomial moment problem for $P(z)$ is either reducible or the
``pull-back'' $q(z)=\tilde q(W(z))W'(z)$ of a solution $\tilde q(z)$ of the polynomial
moment problem for a compositional left factor $\tilde P(z)$ of $P(z)$ and the points $\tilde a=W(a)$ and $\tilde b =W(b).$

If $a,b$ are not critical points of $P(z)$ then $M_{P,a,b}$ always contains elements of form \eqref{compp}. In general case however a more delicate conclusion is true.
Denote by
$P_{a_1}^{-1}(z),$
$P_{a_2}^{-1}(z), ... , P_{a_{d_a}}^{-1}(z)$
(resp. $P_{b_1}^{-1}(z),$ $P_{b_2}^{-1}(z),$ ... ,
$P_{b_{d_b}}^{-1}(z)$)
branches
of $P^{-1}(z)$ in $U$
which map points close
to $P(a)$ (resp. to $P(b)$)
to points close to $a$ (resp. $b$).
In particular, the number $d_a$ (resp. $d_b$)
equals the multiplicity of the point $a$ (resp. $b$)
with respect to $P(z)$. The proposition below was proved in \cite{pry} and by a different method in \cite{pp}.
Below we give a proof following \cite{pp}.

\bp [\cite{pry}, \cite{pp}] \l{brc} If $P(a)=P(b)$ then $M_{P,a,b}$
contains the element
\be \l{e1}
\frac{1}{d_{a}}\sum_{s=1}^{d_{a}}
Q(P_{a_s}^{-1}(z))=\frac{1}{d_{b}}\sum_{s=1}^{d_{b}}
Q(P_{b_s}^{-1}(z)).
\ee
On the other hand, if $P(a)\neq P(b)$ then $M_{P,a,b}$
contains the elements
\be \l{e2}
\frac{1}{d_{a}}\sum_{s=1}^{d_{a}}
Q(P_{a_s}^{-1}(z))=0, \ \ \ \ \ \ \ \ \ \ \frac{1}{d_{b}}\sum_{s=1}^{d_{b}}
Q(P_{b_s}^{-1}(z))=0.
\ee
\ep

\pr Suppose first that $P(a)=P(b).$ Without loss of generality assume that $P(a)= P(b)=c_1$ and consider the relation
$$\phi_1(z)=\sum_{i=1}^nf_{1,i}Q(P^{-1}_{i}(z))=0$$ corresponding to $c_1$.
Let $i,$ $1\leq i \leq n,$ be an index such that $f_{1,i}\neq 0$
and $x$ be a vertex of the star $S_i$ such that $P(x)=c_1.$
It follows from the definition of $\phi_i(z),$ $1\leq i \leq k,$ that if $x\neq a,b$ then there exists an index
$j$ such
that $x$ also is a vertex of the star $S_{j}$ and
$f_{1,j}=-f_{1,i}.$ Furthermore, we have $j=g_1^l(i)$
for some natural number $l$ (see Fig. 2).
\begin{figure}[ht]
\epsfxsize=5truecm
\centerline{\epsffile{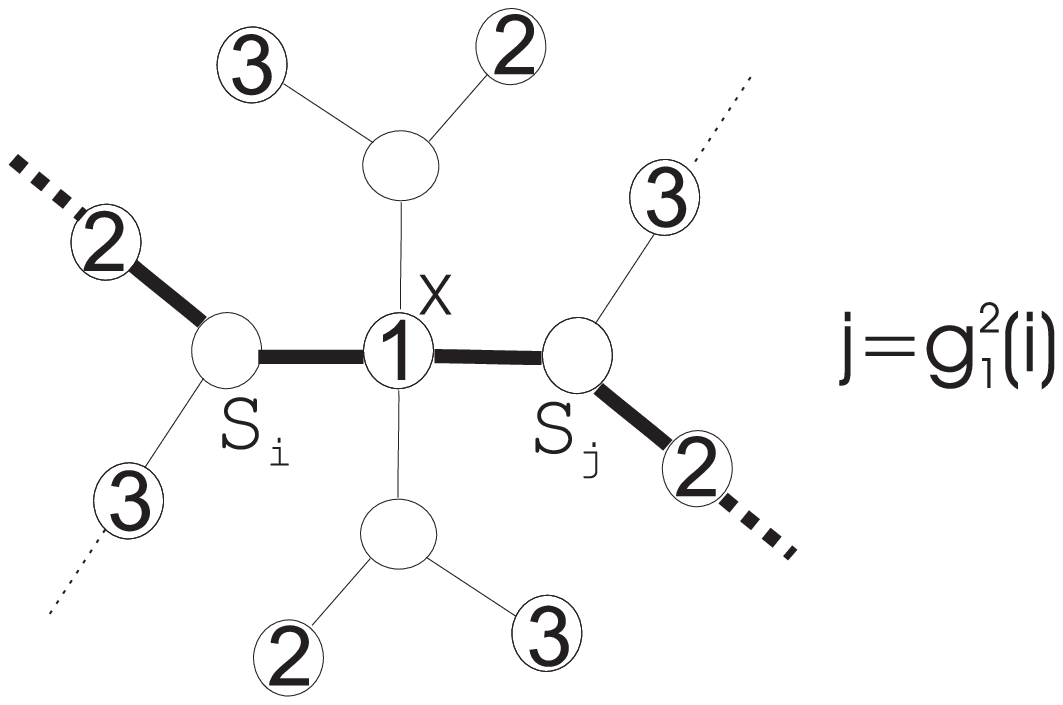}}
\smallskip
\centerline{Figure 2}
\end{figure}
Therefore, $\phi_1(z)$ has the form
$$\phi_1(z)=-Q(P^{-1}_{i_a}(z))+ $$
$$
Q(P^{-1}_{i_1}(z))-Q(P^{-1}_{g_1^{l_1}(i_1)}(z))+
\, ... \, +
Q(P^{-1}_{i_r}(z))-Q(P^{-1}_{g_1^{l_r}(i_r)}(z))$$
$$+Q(P^{-1}_{i_b}(z))=0,$$
where $i_a$ (resp. $i_b$) is an index such that
$a\subset S_{i_a}$
(resp. $b\subset S_{i_b}$), $i_1,i_2,\, ... \, i_r$
are some other indices, and $l_1,l_2,\, ... \, l_r$
are some natural numbers.

For each $s\geq 0$ the equality
$$-Q(P^{-1}_{g^s_1(i_a)}(z))+ $$
$$
Q(P^{-1}_{g^s_1(i_1)}(z))-Q(P^{-1}_{g_1^{l_1+s}(i_1)}(z))+
\, ... \, +
Q(P^{-1}_{g^s_1(i_r)}(z))-Q(P^{-1}_{g_1^{l_r+s}(i_r)}(z))$$
$$+Q(P^{-1}_{g^s_1(i_b)}(z))=0$$ holds by the analytic continuation
of the equality
$\phi_1(z)=0.$ Summing now these equalities from $s=0$ to $s=r-1,$ where $r$ is the order of the element $g_1$ in the group $G_P,$ and taking into account that
for any $i,$ $1\leq i \leq n,$ and any natural number $l$ we have:
$$\sum_{s=0}^{r-1}
Q(P^{-1}_{g^s_1(i)}(z))=\sum_{s=0}^{r-1} Q(P^{-1}_{g_1^{l+s}(i)}(z)),$$
we obtain equality \e{e1}.

In order to prove the proposition in the case when $P(a)\neq P(b)$ it is enough to
examine in a similar way the relations $\phi_1(z)=0$ and $\phi_2(z)=0,$ where $P(a)=c_1,$
$P(b)=c_2.$
\qed

\bc \l{zzaa} Let $P(z),Q(z)\in \C[z]$ and $a,b\in \C,$ $a\neq b.$
Then $\hat H(t)\equiv 0$ implies that \eqref{ux} hold for $q(z)=Q^{\prime}(z).$
\ec
\pr Indeed, if $P(a)=P(b)$ then 
equating the limits of both parts of equality \eqref{e1} as $z$ approaches to $P(a)=P(b)$ 
we see that $Q(a)=Q(b)$. On the other hand, if $P(a)\neq P(b)$ then it follows from equalities \eqref{e2} in a similar way that 
$Q(a)=Q(b)=0$.
In both case it follows from \eqref{ururu} that \eqref{ux} holds.  
\qed

Recall that we assume that the numeration of branches $P_i^{-1}(z),$
$1\leq i \leq n,$ in $U$ is chosen in such a way that the permutation $g_{\infty}\subset G_P$ coincides with the cycle $(1\,2\,...\,n)$.
The proposition below describes the position of branches appearing in Proposition \ref{brc}
with respect to this numeration. More precisely, we describe the mutual position on the unit circle of the sets
$$V(a)=\{ \v_n^{a_1},\v_n^{a_2}, ... , \v^{a_{d_a}}_n \}\ \ \ \ \
{\rm and} \ \ \ \ \ V(b)=\{\v_n^{b_1},\v_n^{b_2}, ... ,\v_n^{b_{d_b}}\},$$
where $\varepsilon_n=exp(2\pi \sqrt{-1}/n)$.

Let us introduce the following
definitions.
Say that two sets of points $X,Y$ on the unit circle $S^1$ are
{\it disjointed} if there exist $s_1, s_2 \in S^1$
such that one of two connected components of $S^1\setminus \{s_1,
s_2\}$ contains all points from $X$ while the other connected component of $S^1\setminus \{s_1,
s_2\}$ contains all points from $Y.$
Say that $X,Y$ are
{\it almost disjointed} if $X\cap Y$ consists of a single point $s_1$
and there exists a point $s_2\in S^1$ such that
one of two connected components of $S^1\setminus \{s_1,
s_2\}$ contains all points from $X\setminus s_1$ while the other connected component of $S^1\setminus \{s_1,
s_2\}$ contains all points from $Y\setminus s_1.$

\bp [\cite{pp}] \l{mono}
The sets $V(a)$ and $V(b)$ are disjointed or almost disjointed.
Furthermore,
if $P(a)= P(b)$ then $V(a)$ and $V(b)$ are disjointed.
\ep
\pr Consider first the case when $P(a)=P(b)=c_1.$
Let $\hat U$ be a simply-connected domain, containing no critical values of $P(z)$,
such that $U\subset \hat U$ and $\infty\in \partial \hat U$. Any branch of $P^{-1}(z)$ in $U$ can be extended analytically to $\hat U$ and we will assume that the numeration of branches of $P^{-1}(z)$ in $\hat U$ is induced by the numeration of
branches of $P^{-1}(z)$ in $U.$
Furthermore,
let $M\subset \hat U$ be a simple curve connecting points $c_1$
and $\infty$ and $\Omega=P^{-1}\{M\}$ be the preimage
of $M$ under the map $P(z)\,:\,\C\P^1\rightarrow
\C\P^1.$
It is convenient to consider $\Omega$ as a bicolored graph
embedded into the Riemann sphere. Namely, we define
black vertices of $\Omega$ as preimages of
$c_1,$ a
unique white vertex of
$\Omega$ as the preimage of $\infty,$ and
edges of $\Omega$ as preimages of $M$ (see Fig. 3).
\begin{figure}[ht]
\medskip
\epsfxsize=10truecm
\centerline{\epsffile{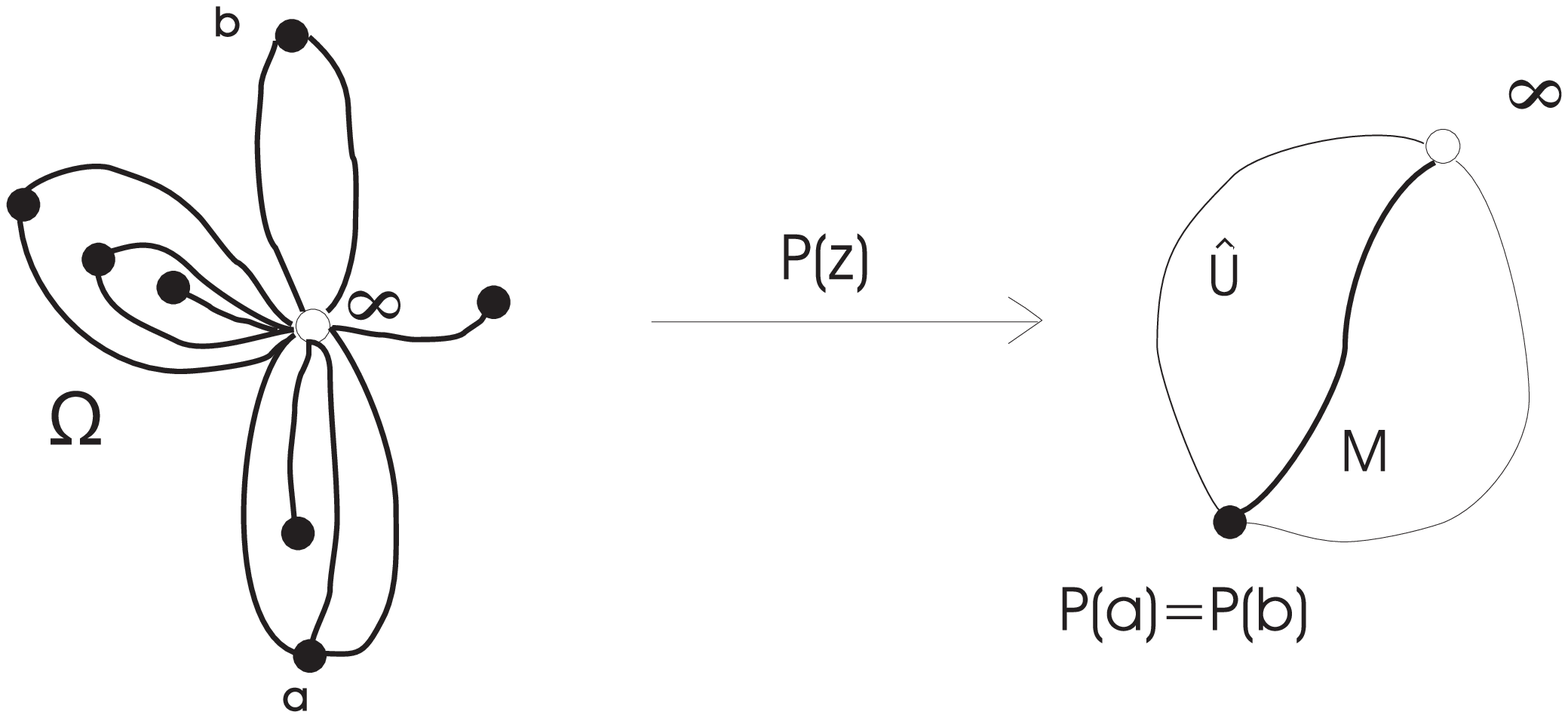}}
\smallskip
\centerline{Figure 3}
\medskip
\end{figure}
The edges of $\Omega$ may be identified with branches of $P^{-1}(z)$ in $\hat U$
as follows:
to the branch $P^{-1}_i(z),$ $1\leq i \leq n,$
corresponds the edge $e_i$
such that $P^{-1}_i(z)$ maps bijectively the interior of $M$ to the interior of $e_i.$
In particular, the ordering of branches of $P^{-1}(z)$ in $\hat U$ induces the ordering of edges of $\Omega.$
Since the multiplicity of the vertex $\infty$ equals $n$ and $\Omega$
has $n$ edges, $\Omega$ is connected.

Let $E_a$
(resp. $E_b$)
be a union of edges of
$\Omega$ which are adjacent to the vertex $a$ (resp. $b$).
It follows from the bijectivity of branches of $P^{-1}(z)$ on the interior of $M$
that if $D$ is a domain from the collection of domains $\C\P^{1}\setminus E_a$ such that $b\in D$, then $D$ contains the
whole set $E_b\setminus \infty.$ Now the proposition follows from the observation that
the cyclic ordering of edges of $\Omega$, induced by
the cyclic ordering of branches of $P^{-1}(z)$ in $\hat U$, coincides
with the cyclic ordering of edges of $\Omega$, induced by the orientation of $\C\P^1$
in a neighborhood of infinity.

In the case when $P(a)\neq P(b)$ the proof is modified as follows.
Take two
simple curves $M_1,$ $M_2\subset \hat U$ connecting the point $\infty$ with the points $P(a)$ and $P(b)$ correspondingly
and consider the preimage $P^{-1}\{M_1\cup M_2\}$
as a graph
$\Omega$ embedded into the Riemann sphere. The vertices of $\Omega$
fall into three sets: the first one consists of a unique vertex
which is the preimage of $\infty,$ the second one consists of vertices
which are preimages of $P(a),$ and the third one consists of vertices
which are preimages of $P(b).$ Similarly, the edges of $\Omega$
fall into two sets: the first one consists of edges
which are preimages of $M_1$ and the second one consists of edges
which are preimages of $M_2$
(see Fig. 4).
\begin{figure}[ht]
\medskip
\epsfxsize=9truecm
\centerline{\epsffile{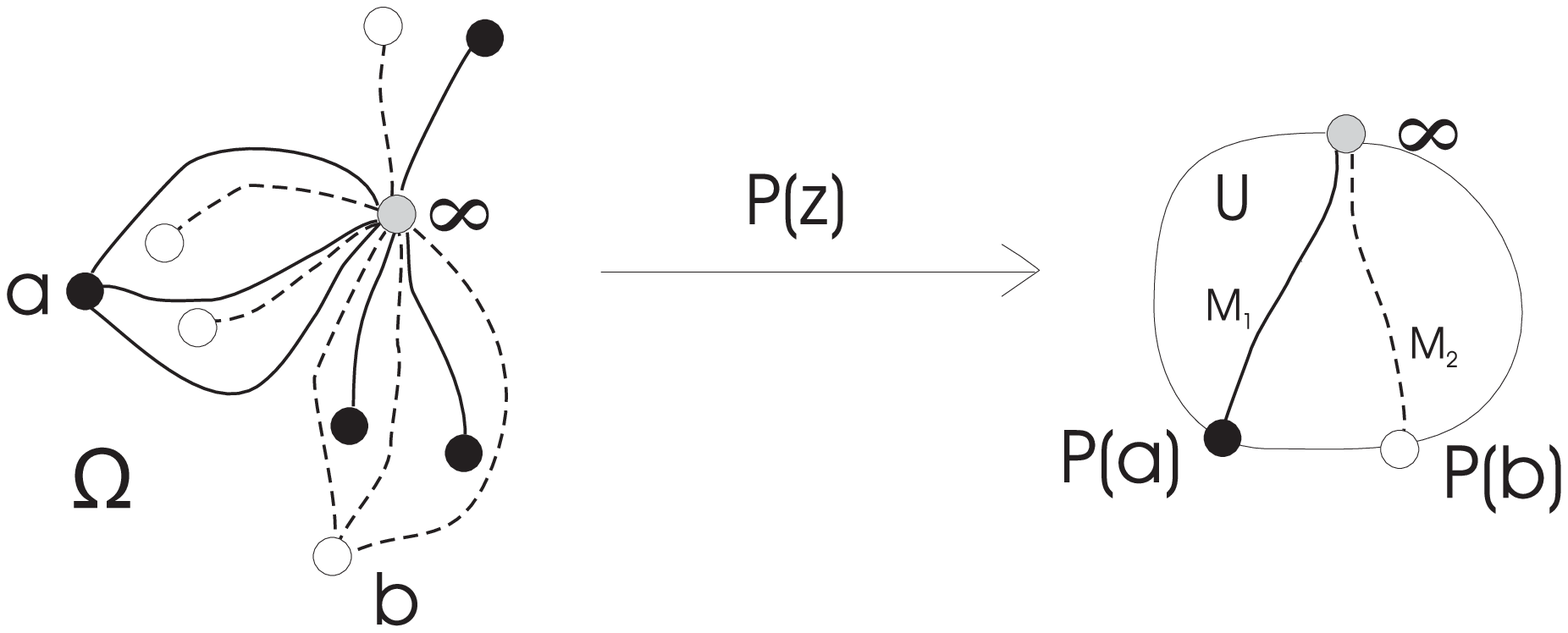}}
\smallskip
\centerline{Figure 4}
\medskip
\end{figure}

Each of two sets of edges of $\Omega$ may be identified with branches of $P^{-1}(z)$ in $\hat U$
as follows:
to the branch $P^{-1}_i(z),$ $1\leq i \leq n,$
corresponds the edge $e^1_i$ from the first set (resp. the edge $e^2_i$ from the second set)
such that $P^{-1}_i(z)$ maps bijectively the interior of $M_1$ (resp. of $M_2$) to the interior of $e^1_i$
(resp. of $e^2_i$).
The ordering of branches of $P^{-1}(z)$ in $\hat U$ induces the ordering of edges of $\Omega$
in each of two sets. Clearly, this
ordering coincides with the natural ordering induced by the
orientation of $\C\P^1.$ Furthermore, it is easy to see that
when going round infinity in the counterclockwise direction the edge
$e^1_i,$ $1\leq i \leq n,$ is
followed by the edge $e^2_i.$

Let $E_a^1$
(resp. $E_b^2$)
be a union of edges from the first (resp.
the second) set
$\Omega$ which are adjacent to the vertex $a$ (resp. $b$). The bijectivity of branches of $P^{-1}(z)$ on the interior of $M_1$ and $M_2$
implies that if $D$ is a domain from the collection of domains $\C\P^{1}\setminus E_a^1$ such that $b\in D$, then $D$ contains the
whole set $E_b^2\setminus \infty.$ 
Taking into account that
for any $k,$ $1\leq i \leq n,$
the edge $e^1_i$ is
followed by $e^2_i,$ this implies
that $V(a)$ and $V(b)$ are disjointed or almost disjointed.
\qed

\noindent{\bf Remark.} Since $Q(P^{-1}_{i}(z))$, $1\leq i \leq n,$ are branches of an algebraic function,
relations \eqref{sss} are examples of linear relations between roots of an algebraic equation over the field $\C(z).$ 
A general algebraic approach to such relations, 
over an arbitrary field, was developed in the papers \cite{gi1}, \cite{gi2}. In particular,
it follows from Theorem 1 of \cite{gi2}
that a necessary and sufficient condition for the existence of at least one solution $Q(z)$ of \eqref{sss},
such that the functions $Q(P^{-1}_{i}(z))$, $1\leq i \leq n,$ are distinct between themselves, is that the subspace $M_{P,a,b}$ does not contain elements of form \eqref{gi}. An equivalent form of this condition is that the subspace $M_{P,a,b}$ does not contain any of subspaces $V_d^{\perp},$ $d\in D(G_P),$ which are defined below.
Notice however that the method of \cite{gi2} does not provide any information about the
description or the actual finding of these solutions.

\section{\l{cycy} Permutation matrix representations of groups containing a full cycle}
\subsection{Invariant subspaces and the centralizer ring}

The construction of $M_{P,a,b}$ implies that $M_{P,a,b}$ is an invariant subspace
of $\Q^n$ with respect to the so called {\it permutation matrix representation}
of the group $G_P$ on $\Q^n.$
By definition, the
permutation matrix representation of a transitive permutation group $H\subseteq S_n$
on $\Q^n$ is a homomorphism $R_H:\, H\rightarrow GL_n(\Q)$
which associates to
$h\in H$ a permutation matrix $R_H(h)\in GL_n(\Q)$ the elements $r_{i,j},$ $1\leq i,j\leq n,$ of  
which satisfy $r_{i,j}=1$ if $i=j^h$ and $r_{i,j}=0$ otherwise. In other words,
$$\begin{pmatrix}x_1\\x_2\\ \vdots \\x_n\end{pmatrix}=R_H(h)
\begin{pmatrix}x_{1^h}\\ x_{2^h}\\ \vdots \\ x_{n^h}\end{pmatrix}.$$
Note that $\Q^n$ admits a $R_H$-invariant scalar product $(x, y):=\sum_{i=1}^n x_iy_i$.

The goal of this section is to provide a full description of the invariant subspaces of $\Q^{n}$
with respect to the permutation matrix representation of $G_P$. More general, we classify all invariant subspaces of $\Q^{n}$
with respect to the permutation matrix representation of an arbitrary  group
$G\subseteq S_n$ containing the cycle $(1\,2\,...\,n)$. {\it In the following $G$ will always denote such a group.}

Recall that a subset $B$ of $X=\{1,2,\dots ,n\}$
is called a {\it block} of a transitive permutation group $H\subseteq S_n$ if for each $h\in H$ the set
$B^h$ is either disjoint or is equal to $B$ (see e.g. \cite{wi}). For a block $B$ the set
${\mathcal B}:=\{B^h\,|\,h\in H\}$ forms a partition of $X$ into a disjoint union
of blocks of equal cardinality which is called an {\it imprimitivity system} of $H$.
Each permutation group $H\subseteq S_n$ has two {\it trivial} imprimitivity systems:
one formed by singletons and another formed by the whole $X$. A permutation group
is called {\it primitive} if it has only trivial imprimitivity systems. Otherwise it
is called {\it imprimitive}.

For each $d\,|\,n$ we denote by $V_d$ the subspace of $\Q^{n}$
consisting of $d$-periodic vectors. The fact that the group $G$ contains the cycle $(1,...,n)$
implies easily the following statement.

\bl\label{p_1} Any imprimitivity system for $G$ coincides with the residue classes modulo $d$
for some $d\,|\,n$. Furthermore, for given $d$ such classes
form an imprimitivity system for $G$
if and only if the subspace $V_d$ is $G$-invariant.\qed
\el

Denote by $D(G)$ the set of all divisors of $n$ for which $V_d$ is $G$-invariant. Clearly, $1,n\in D(G)$. Notice that $D(G)$ is a
lattice with respect to the ope\-rations $\land, \lor,$ where
$d\land f := \gcd(d,f)$ and $d\lor f := {\rm lcm}(d,f)$. Indeed, for an element $x\in X$
the intersection of two blocks containing $x$ and corresponding to $d,f\in D(G)$ is a block which
corresponds to $d\lor f$. On the other hand,
the intersection of two invariant subspaces $V_d,V_f$ is
an invariant subspace which is equal to $V_{d\land f}.$

We say that $d\in D(G)$ {\it covers} $f\in D(G)$
if $f\,|\,d,$ $f<d,$ and there is no $x\in D(G)$ such that $f<x<d$ and $f\vert x,$ $x\vert d$.
Now we are ready to formulate the main result of this section.

\begin{thm}\label{p_main} Each $R_G$-irreducible subspace of $\Q^{n}$ has the form
$$U_d:=V_d\cap \left(V_{f_1}^\perp\cap ...\cap V_{f_\ell}^\perp\right),$$ where $d\in D(G)$ and
$f_1,...,f_\ell$ is a complete set of elements of $D(G)$ covered by $d$. The subspaces $U_d$ are mutually orthogonal and
every $R_G$-invariant subspace of $\Q^{n}$ is a direct sum of some $U_d$ as above.
\end{thm}

The proof of this theorem splits into 
several steps and is given below.
We start from recalling some basic 
facts of the representations theory
which we will use afterward (see e.g. \cite{kir}).

First, any representation $T_H:\, H\rightarrow GL_n(k)$ of a finite group $H$ over a field $k$ of characteristic not dividing
$\vert H\vert$ 
is completely reducible,
that is $k^n$ is a direct sum of $T_H$-invariant irreducible  subspaces (Maschke's theorem).
Furthermore, irreducible subspaces
of a completely reducible representation $T_H:\, H\rightarrow GL_n(k)$
are in one-to-one correspondence with minimal idempotents of the {\it centralizer ring} $V_{k}(T_H)$.
Recall that $V_{k}(T_H)$ consists of all matrices $A\in M_n(k)$ which commute
with every $T_H(h), h\in H.$ Furthermore, a matrix $E$ is called
an idempotent if $E\neq 0$ and $E^2 = E$.
Two idempotents $E,F$ are called {\it orthogonal} if $EF = FE = 0$. Finally, an idempotent
$E\in V_{k}(T_H)$ is called {\it minimal} if it can not be presented as
a sum of two orthogonal idempotents from $V_{k}(T_H)$. Under this notation the correspondence above is obtained as follows: to a minimal idempotent $E\in V_{k}(T_H)$ corresponds an irreducible subspace $V=\im\{E\}.$

In general, the decomposition of $k^n$
into a sum of $T_H$-invariant irreducible subspaces is not uniquely defined. Nevertheless,
if \be \l{dec}
k^n=V_1^{\oplus a_1}\oplus\dots\oplus V_r^{\oplus a_r}\ee is a decomposition such that
$V_i,$ $1\leq i \leq r,$ are pairwise non-isomorphic 
$T_H$-invariant irreducible subspaces of $k^n$,
then the
subspaces $V_i^{\oplus a_i},$ $1\leq i \leq r,$ are defined uniquely. They
correspond to the
minimal idempotents of the {\it center} $C(V_{k}(T_H))$ of 
the centralizer ring $V_{k}(T_H).$ Furthermore,
$V_{k}(T_H)$ is commutative if and only if
$a_i = 1$ for all $i,$ $1\leq i\leq r$.
Notice that if $V_{k}(T_H)$ is commutative and the space $k^n$ admits a $T_H$-invariant scalar product
then all $T_H$-invariant irreducible subspaces of $k^n$ are mutually orthogonal. Indeed, for any representation $T_H:\, H\rightarrow GL_n(k)$, which admits an invariant scalar product, $k^n$ can be decomposed into a sum of 
$T_H$-invariant irreducible subspaces
\be \l{dec1}
k^n=V_1\oplus\dots\oplus V_r\ee with mutually orthogonal $V_i$. On the other hand, if $V_{k}(T_H)$ is commutative then a decomposition of $T_H$
into a sum of $T_H$-invariant irreducible subspaces is uniquely defined and therefore coincides with \eqref{dec1}.

For the permutation matrix representation $R_H:\, H\rightarrow GL_n(k)$ of a transitive permutation group $H\subseteq S_n$ instead of the notation $V_{k}(R_H)$ we will use simply the symbol $V_{k}(H).$
Below we will show (Proposition \ref{corre}) that for any group $G$ as above the ring
$V_{\Q}(G)$ is isomorphic to a subring of the group algebra of a cyclic group
and hence is
commutative. Therefore, the above remarks imply the following statement.

\bp \label{idempotent}
An $R_G$-invariant subspace $W\subset{\Q}^n$ is irreducible
if and only if there exists a minimal idempotent
$E\in V_{\Q}(G)$ such that $\im\{E\} = W$. $R_G$-invariant irreducible subspaces of $\Q^n$ are mutually orthogonal and every $R_G$-invariant subspace is a
direct sum of some $W$ as above. \qed
\ep

For each transitive permutation group $H\subseteq S_n$ we can construct
some special basis of $V_\C(H)$ via orbits of the stabilizer
$H_1$ of the point 1 as follows. To each orbit $\Delta$ of $H_1$ associate a matrix $A^{\Delta}$,
where $A_{i,j}^{\Delta}=1$ if there exist $h\in H$, $\delta \in \Delta$ such that $1^h=j,$ $\delta^h=i,$ and
$A_{i,j}^{\Delta}=0$ otherwise. In particular, for the first column
of $A^{\Delta}$ the equality $A_{i,1}^{\Delta}=1$ holds if and only if $i\in \Delta.$
It turns out that the matrices $A^{\Delta}$ form a basis of $V_\C(H)$ (\cite{wi}, Theorem 28.4). Furthermore,
since by construction the matrices $A^{\Delta}$ are contained in
$M_n(\Q)$ they form a basis of $V_\Q(H)$.
We summarize the properties of $A^{\Delta}$ in the proposition below (see \cite{wi}, \S 28).

\bp \l{alg} The matrices $A^{\Delta}$ satisfy the following conditions:
\begin{itemize}

\item[(1)] $A^{\Delta}$ form a basis of the algebra $V_\Q(H)$ as of a $\Q$-module;

\item[(2)] If $\Delta_1\neq \Delta_2$ then the ones of $A^{\Delta_1}$ and $A^{\Delta_2}$
do not occur in the same place. On the other hand, $\sum_{\Delta}A^{\Delta}$ is a matrix all the entries of
which are ones;

\item[(3)] For each orbit $\Delta$ there exists an orbit $\Gamma$ such that
$(A^{\Delta })^T=A^{\Gamma}.$ \qed 

\end{itemize}

\ep
Notice that the property (3) implies that for the first row of
$A^{\Delta}$ the equality $A_{1,j}^{\Delta}=1$ holds if and only if $j\in \Gamma$.
Furthermore, it is easy to see that the mapping $\Delta\rightarrow \Gamma$ defines an involution
on the set of orbits of $H_1.$

\pagebreak

\subsection{Schur rings}

\subsubsection{Isomorphism between $S_{\Q}(G)$ and $V_\Q(G)$}
In order to construct the minimal idempotents of $V_{\Q}(G)$ we will use
so called {\it Schur rings} introduced by Schur in his classical paper \cite{schu} for the investigation of
permutation groups
containing a regular subgroup $C$.
Since in this paper $C$ always will be a cyclic group, in the following we will restrict our attention to
this case only (see \cite{wi} for the account of the Schur method in the general case).

The idea of the Schur approach can be described as follows.
If $G\subseteq S_n$ contains the cycle $c:=(1\,2\, ...\,n)$
then elements of the set $\left\{1,2,\dots,n\right\}$ can be identified with elements of the cyclic group $C$ generated
by $c$ as follows: to the element $i$ corresponds the element of $C$ which transforms $1$ to $i$.
Therefore, we can consider $G$ as a permutation group acting on its subgroup $C.$
After such an identification we can ``multiply'' elements of the set $\left\{1,2,\dots,n\right\}$ and
this multiplication agrees with the action of $G$ in the following sense: if $h, g \in C$ then $h^g=hg$.
Furthermore, identifying any two subsets of $\left\{1,2,\dots,n\right\}$
with the corresponding elements of the group algebra $\Q[C]$ we can define their ``pro\-duct''
as the product of these elements in $\Q[C]$. The remarkable result of Schur is that under
such a multiplication the orbits of the stabilizer $G_1$
form a basis of some subalgebra
of $\Q[C]$.
To make this statement precise let us introduce the following definition.

For $T\subseteq C$
denote by $T^{\left(-1\right)}$ the set of elements of $C$ inverse to the elements of $T$
and by $\und{T}$ the formal sum $\sum_{h\in T} h$. The elements of $\Q[C]$ of the form $\und{T}$ for
some $T\subseteq C$
are called {\it simple quantities} (\cite{wi}).
\begin{definition}\label{s-ring}
A subalgebra $\cA$ of the group algebra $\Q [C]$ is called a
{\it Schur ring} or an {\it S-ring} over $C$
if it satisfies the following axioms:

\begin{itemize}

\item[(S1)] $\cA$ as a $\Q$-module has a basis consisting
         of simple quantities
        $\und{T_0},\dots,\und{T_d}$, where $T_0=\{e\}$,

\item[(S2)] $T_i\cap T_j =\emptyset$ for $i\neq j$ and
 $\bigcup_{j=0}^d T_j=C$,

\item[(S3)] For each $i\in\{0,1,\dots,d\}$ there exists
$i'\in\{0,1,\dots,d\}$ such that \linebreak
$T_{i'} =T_i{}^{(-1)}$.

\end{itemize}
\end{definition}

It is easy to see that the basis $\und{T_0},\dots,\und{T_d}$ satisfying (S1) and (S2) is 
unique. Such a basis is called the {\it standard basis} of $\cA$.
The number
$d+1$ is called the {\it rank} of $\cA$. The sets $T_i$, $0\leq
i\leq d$, are called the {\it basic sets} of $\cA$. Finally, the
notation $\cA=\langle  \und{T_0},\dots,\und{T_d} \rangle$ is
used if $\cA$ is an S-ring over $C$ whose basic sets are
$T_0,\dots,T_d$.
We also write $\Bsets{\cA}$ for the set
$\{T_0,\dots,T_d\}$.
Notice that if $\tilde \cA$ is an $S$-ring which is a subring of
$\cA$ then its basic sets are some unions of basic sets of $\cA$. There
are two {\it trivial} S-rings, namely
$\langle \und{e},\und{C\setminus\{e\}}\rangle$ and $\Q[C]$.

\bp \l{corre} To any group $G$ corresponds a Schur ring $S_{\Q}(G)$ the basic sets of which
are the orbits of the stabilizer $G_1$. Moreover,
$S_{\Q}(G)$ and $V_\Q(G)$ are isomorphic as $\Q$-algebras.

\ep

The Proposition \ref{corre} is a particular case of Theorem 28.8 in \cite{wi}. It
implies in particular that in order to describe
the minimal idempotents of $V_{\Q}(G)$ it is enough to describe the ones of $S_{\Q}(G)$. Since however
for this purpose
an explicit construction of the isomorphism between $S_{\Q}(G)$ and $V_\Q(G)$
is needed, we give below a short proof of Proposition \ref{corre} which
is based on Proposition \ref{alg}

\vskip 0.2cm
\noindent{\it Proof of Proposition \ref{corre}.} First of all observe that
since $G$ contains $c$
each matrix $M\in V_\Q(G)$ is necessarily a {\it circulant} that is
each row vector of $M$
is cyclically shifted for one element to the right relative to the preceding row vector, in other words
\be \l{cirrr} M_{i,j} = M_{1,j-i+1 \, \mod n}.\ee

Define now a mapping $\psi: V_\Q(G)\rightarrow \Q[C]$
by the formula $$\psi(M):=\sum_{j=1}^{n} M_{1,j}c^{j-1}$$ and show
that $\psi$ is an algebra monomorphism.
Indeed, for any $M,N\in V_\Q(G)$ we have:
$$
\psi(MN) = \sum_{\ell=1}^{n} (MN)_{1,\ell}c^{\ell-1} = \sum_{\ell=1}^{n} \sum_{i=1}^{n}
M_{1,i}N_{i,\ell} c^{\ell-1}=$$
$$=
\sum_{\ell=1}^{n} \sum_{i=1}^{n} M_{1,i}N_{1,\ell-i+1} c^{\ell-1}=
\sum_{i=1}^{n} \sum_{j=1}^{n}M_{1,i}N_{1,j}c^{i+j-2} =$$ $$ =\left(\sum_{i=1}^{n} M_{1,i}c^{i-1}\right)
\left(\sum_{j=1}^{n} N_{1,j}c^{j-1}\right) =
\psi(M)\psi(N).
$$
Thus $\psi$ is an algebra homomorphism. Furthermore, $\psi$ is injective since
any matrix $M\in V_\Q(G)$ is defined by its first row in view of \eqref{cirrr}.

Clearly, the image of $V_\Q(G)$ is a subalgebra $S_{\Q}(G)$ of $\Q[C]$. Furthermore, by construction, the basis of this subalgebra consists of the orbits of the stabilizer
$G_1.$ The properties S1, S2 of $S_{\Q}(G)$ are obvious. Finally, since
any matrix from $V_\Q(G)$ is a circulant, it follows from the third part of Proposition \ref{alg} that $\Delta^{(-1)}=\Gamma$. \qed

For $d$ dividing
$n$ denote by $C_d$ a unique subgroup of $C$ of order $d$. For a Schur ring
$\cA$ denote by $D(\cA)$ a set consisting of all divisors of $n$ for which $\und{C_d}\in\cA $.

\bl\label{p_D}
$d\in D(G)\iff n/d\in D(S_{\Q}(G))$.
\el
\pr Let $d\in D(G)$. Then $C_{n/d}$ under the identification of the set $\{1,2, \dots, n\}$ with $C$
corresponds to
the set $X=\{1,d+1,2d+1, \dots, n-d+1\}$ and therefore
is a block of $G$ containing 1. This implies that $C_{n/d}$ is a union of some $G_1$-orbits,
say $T_0,...,T_\ell$. Hence
$\und{C_{n/d}}= \und{T_0}+\und{T_1}+\dots+\und{T_\ell}$
and therefore $\und{C_{n/d}}\in S_{\Q}(G)$.

Let now $n/d\in D(S_{\Q}(G))$. Then $\psi^{-1}(\und{C_{n/d}})\in V_\Q(G)$.
It follows from the definition of $\psi$ that $\psi^{-1}(\und{C_{n/d}})$
is a circulant matrix $M$ such that $M_{1,i}=1$ if $i\in X$ and
$0$ otherwise. Since
$M \in V_\Q(G)$ the subspace ${\rm Im}(M)$
is $G$-invariant. On the other hand, it is easy to see that
$ {\rm Im}(M)= V_d$. Therefore, $d\in D(G)$ by Lemma \ref{p_1}. \qed

\subsubsection{Rational $S$-rings}
The automorphism group
of $C$ is isomorphic to the multiplicative group $\Z_n^*$.
Namely, to the element $m\in\Z_n^*$ corresponds the automorphism $g\mapsto g^m,g\in C.$
Extending this action onto $\Q[C]$
by linearity we obtain an action of $\Z_n^*$ on the group algebra $\Q[C]$:
$$\alpha=\sum_{g\in C} \alpha_g g\ \ \longrightarrow \ \alpha^{(m)}:=\sum_{g\in C} \alpha_g g^m.$$ An element $\alpha\in\Q[C]$ is called {\it rational} if $\alpha=\alpha^{(m)}$ for any $m \in \Z_n^*$.
Note that the mappings $\alpha\mapsto\alpha^{(m)},$ $m \in \Z_n^*$,
are automorphisms of $\Q[C]$. Moreover, these mappings
are automorphisms of any S-ring $\cA$ over $C$
(see \cite{wi}, Theorem 23.9). In particular,
for each $m\in\Z_n^*$ and $T\subseteq C$ we have
$$
T\in\Bsets{\cA}\iff T^{(m)}\in\Bsets{\cA},
$$
where for a subset $T\subset C$ by $T^{(m)}$ is denoted the set of $m$-th powers of
elements of $T$.

Recall that the set of all irreducible complex representations of $C$ consists of $n$ one-dimensional representations
(characters)
$\chi_0,...,\chi_{n-1}$ where
$$\chi_\ell(c^j):=e^{2\pi \sqrt{-1} \ell j/n}, \ \ 0\leq j,\ell\leq n-1.$$
We will keep the same notation for the extensions of $\chi_0,...,\chi_{n-1}$ by linearity
on $\Q[C]$.
The rational elements of an $S$-ring $\cA$ admit the following characterization.

\bl\l{re}
An element $\alpha\in\Q[C]$ is rational if and only if
$\chi_l(\alpha)\in\Q$ for all $l,$ $0 \leq l \leq n-1.$
\el

\pr For an element $\alpha=\sum_{j=1}^n h_j c^j$ of $\Q[C]$ the condition
that $\chi_l(\alpha)\in\Q$ for all $l,$ $0 \leq l \leq n-1,$ is equivalent to the condition that
$\chi_l(\alpha),$ $0\leq l \leq n-1,$
is invariant with respect to the action of the Galois group $\Gamma$ of the extension
$(\Q(e^{2\pi \sqrt{-1} /n}):\Q)$. The group $\Gamma$ is isomorphic to $\Z_n^*$.
Namely, to the element $m\in \Z_n^*$
corresponds the element $\sigma_m\in \Gamma$ which transforms $e^{2\pi \sqrt{-1} /n}$ to $e^{2\pi \sqrt{-1} m /n}$.
We have:
$$ \sigma_m(\chi_\ell(\alpha)) = \sigma_m(\chi_\ell(\sum_{j=1}^n h_j c^j))
=\sigma_m(\sum_{j=1}^n h_j e^{2\pi \sqrt{-1} \ell j/n})
=
$$$$=\sum_{j=1}^n h_j e^{2\pi \sqrt{-1} m\ell j/n}=
\chi_{\ell}(\sum_{j=1}^n h_j c^{mj}) = \chi_\ell(\alpha^{(m)}).$$
Therefore, for $\ell,$ $0\leq \ell \leq n-1,$ and
$m\in\Z_n^*$ the equality
$\sigma_m(\chi_\ell(\alpha))=\chi_\ell(\alpha)$ is equivalent to the equality
$\chi_\ell(\alpha^{(m)})=\chi_\ell(\alpha)$. Since
for $\alpha, \beta \in \Q[C]$ the equality  $\chi_{\ell}(\alpha)=\chi_{\ell}(\beta)$ holds
for all $\ell,$ $0\leq \ell \leq n-1,$
if and only if $\alpha=\beta,$ we conclude that $\chi_\ell(\alpha)\in\Q$
for all $\ell,$ $0 \leq \ell \leq n-1,$ if and only if $\alpha$ is rational.\qed

An $S$-ring $\cA$ is called {\it rational} if all its elements are rational.
Clearly, $\cA$ is rational if and only if $T^{(m)}=T$ for all $T\in\Bsets{\cA}$ and $m\in \Z_n^*$.
Any rational $S$-ring is
a subring of some universal rational $S$-ring $W.$
To construct $W$ observe that the orbits of the action of $\Z_n^*$ on
$C$ are parametrized by the divisors of $n$ as follows: an orbit $O_m,$ $m\vert n,$ consists of all generators of
the group $C_m$. It turns out that the vector space spanned by $\und{O_m},$ $m\vert n$,
is a rational $S$-ring $W$ (\cite{schu}). Furthermore, any rational $S$-ring $\cA$ is
a subring of $W.$ Indeed, since any element of the standard basis of a rational $S$-ring $\cA$ is
invariant with respect to the action of $\Z_n^*,$ such an element is
a union of some $O_m,$ $m\vert n.$ Therefore,
$\cA$ is
a subring of $W.$

Denote by $D_n$
the lattice of all divisors of $n$ with respect to the
operations $\land, \lor.$ The statement below describes the rational S-rings.

\bp \l{mu} (\cite{mu})\label{p_rational}
An S-ring $\cA$ over $C$ is rational if and only if there exists a sublattice $D$ of $D_n$
with $1,n\in D$ such that $\und{C}_d,\, d\in D$, is a basis
of $\cA$. \qed
\ep

Notice that the basis $\und{C_d},$ $d\in D$, is not a standard
basis of $\cA$ in the sense of definition \ref{s-ring}.

To any $S$-ring $\cA$ one can associate a rational $S$-ring
$\Trace{\cA}$, called the {\it rational closure} of $\cA$,
which is constructed as follows.
Introduce an equivalence relation on $\Bsets{\cA}$ setting
$S\sim T$ if there exists
$m\in\Z_n^*$ such that $S=T^{(m)}$. For $T\in\Bsets{\cA}$ set $$\Trace{T}:=\bigcup\{T^{(m)}\,|\,m\in\Z_n^*\}$$ and denote by $\Trace{\cA}$ the $\Q$-module spanned by $\und{\Trace{T}},$ $T\in\Bsets{\cA}$.

\bp \l{mmp} (\cite{schu}) $\Trace{\cA}$ is an S-ring consisting of all rational elements of $\cA.$
\ep

The Proposition \ref{mu} allows us to describe a rational closure of an arbitrary S-ring.

\bp \label{p_r_closure2} Let $\cA$ be an S-ring over $C$. Then $\und{C}_d\,,d\in D(\cA)$, is a basis
of $\Trace{\cA}$.
\ep
\begin{proof} By Proposition ~\ref{p_rational} $\Trace{\cA}$
is spanned by vectors $\und{C}_d,\, d\in D$,
for a certain sublattice $D$ of $D_n$. It remains to prove that $D=D(\cA)$. The inclusion
$D\subseteq D(\cA)$ follows from the following line
$$
d\in D \implies \und{C_d}\in\Trace{\cA}\subseteq \cA\implies \und{C_d}\in \cA\implies d\in D(\cA).
$$

Conversely, pick an arbitrary $f\in D(\cA)$. Then $\und{C_f}\in\cA$. Furthermore, since $$\und{C_f}=\sum_{t\in D_f} \und{O_t}\ ,$$ the element
$\und{C_f}$ is rational and therefore
$\und{C_f}\in\Trace{\cA}$. This means that $\und{C_f}$ is
a linear combination of $\und{C_d},\, d\in D$.
Therefore, in order to prove that $\und{C_f} = \und{C_d}$ for suitable $d\in D$ it is enough to show that the simple quantities $\und{C_d},\, d\in D_n,$ are linearly independent.

In order to prove the last statement assume that \be \l{asd} \sum_d l_d\,\und{C_d}=0\ee and let $M$ be a maximal number $d$
for which $l_d\neq 0$. Clearly, any element $u$ of $C$ which generates $C_M$ can not be an element of $C_d$ for $d<M.$ But then $u$ appears in the left part of equality \eqref{asd} only once with the coefficient $l_d\neq 0$. This is a contradiction and therefore
$\und{C_d},\, d\in D_n,$ are linearly independent.
\end{proof}

\subsection{Proof of Theorem \ref{p_main}}
Similarly to the definition given above for the elements of $D(G)$
say that for an $S$-ring $\cA$ an element $d\in D(\cA)$ covers an element $f\in D(\cA)$
if $f\,|\,d,$ $f<d,$ and there is no $x\in D(\cA)$ such that $f<x<d$ and $f\vert x,$ $x\vert d$. 

Set $$\sig_d:=\frac{1}{d}\und{C_d}, \ \ \ d\in D(\cA).$$ It follows from  
\begin{equation}\label{eq_idem}
\sig_f\sig_d = \sig_d\sig_f=\sig_{f\lor d}
\end{equation} that $\sig_d, d\in D(\cA),$ are idempotents of the
algebra $\cA$. Nevertheless, they are not pairwise orthogonal.

\bp \label{p_idempotents} An element of an S-ring $\cA$ over $C$ is a
minimal idempotent of $\cA$ if and only if it has the form
\be \label{minimal_idem}
\eps_d = \sig_d\prod_{i=1}^\ell(1-\sig_{f_i}),
\ee
where $d\in D(\cA)$ and $f_1,...,f_\ell$ is a complete set of elements of $D(\cA)$
covering $d$.
\ep
\begin{proof}
Let us show first that $\eps_d,d\in D(\cA),$ are pairwise orthogonal idempotents.
Since each $\sig_d,$ $d\in D_n$, is an idempotent, we have:
$$
\eps_d^2 = \sig_d^2\prod_{i=1}^\ell(1-\sig_{f_i})^2 =\sig_d\prod_{i=1}^\ell(1-2\sig_{f_i}+\sig_{f_i}^2)=
\sig_d\prod_{i=1}^\ell(1-\sig_{f_i}) = \eps_d.
$$
Therefore, in order to show that $\eps_d$ is an idempotent we only must check that $\eps_d\neq 0.$
In view of \eqref{eq_idem},
after opening the brackets in \eqref{minimal_idem} we obtain
a linear combination of $\sigma_f$ in which $\sig_d$ appears
with the coefficient one. Since $\sigma_d,$ $d\in D_n$, are linearly independent this implies that
$\eps_d\neq 0.$

Let us check now the orthogonality. Take two distinct $m,d\in D(\cA),$
where it is assumed that $d<m,$ and
consider the product $\eps_d\eps_m$. Let $f_1,...,f_\ell$ and $n_1,...,n_k$
be complete sets of elements of $D(\cA)$ which cover $d$ and $m$ respectively.
By \eqref{eq_idem} we have:
$$
\eps_d\eps_m = \sig_d\prod_{i=1}^\ell(1-\sig_{f_i})\cdot\sig_m \prod_{j=1}^k(1-\sig_{n_j})
=\sig_d \sig_m\prod_{i=1,j=1}^{i=\ell,j=k}(1-\sig_{f_i})(1-\sig_{n_j}) =$$
\be \l{lop} =\sig_{d\lor m}\prod_{i=1,j=1}^{i=\ell,j=k}(1-\sig_{f_i})(1-\sig_{n_j}) \ee
Since $d\,|\,d\lor m$ and $d < d\lor m$, there exists an element $f_i\in D(\cA)$ which
covers $d$ and divides $d\lor m$. For such an element
$(1 -\sig_{f_i})\sig_{d\lor m} = 0$ and
this implies the vanishing of the right-hand side of \eqref{lop}.

Since the idempotents $\eps_d,d\in D(\cA),$ are pairwise orthogonal
they are linearly independent elements of $\cA$. Furthermore, since Proposition \ref{p_r_closure2} implies that 
$\eps_d \in \Trace{\cA}$ for any $d\in D(\cA)$ and
\be \l{dim} \dim(\Trace{\cA}) = |D(\cA)|,\ee the idempotents $\eps_d,d\in D(\cA),$ form a basis of $\Trace{\cA}$ which consists of
pairwise orthogonal idempotents.
This implies that any minimal idempotent $\eps$ of $\Trace{\cA}$ coincides with some $\eps_d,d\in D(\cA).$ Indeed, since $\eps_d,d\in D(\cA),$ form a basis of $\Trace{\cA}$ there exist numbers $a_d, d\in D(\cA)$, such that $\eps=\sum_{d\in D(\cA)} a_d\eps_d.$ 
Furthermore, since $\eps$ is an idempotent, for any $d\in D(\cA)$ the coefficient
$a_d$ equals either $1$ or $0.$ Therefore, if $\eps$ is minimal then $\eps=\eps_d$ for some $d\in D(\cA).$

Finally, observe that the sets of minimal idempotents of $\Trace{\cA}$ and $\cA$ coincide. Indeed,
if $\eps$ is any idempotent of $\cA$ then $\eps^2 = \eps$ implies that $\chi_i(\eps)\in\{0,1\}$ for all $i,$ $0\leq i \leq n-1.$ Therefore, by Proposition \ref{mmp},
$\eps\in \Trace{\cA}.$ Furthermore, if $\eps$ is minimal in $\cA$ then obviously it is also minimal
in $\Trace{\cA}$. On the other hand, any minimal idempotent of  $\Trace{\cA}$ remains a minimal idempotent in $\cA$ since all idempotents of $\cA$ are contained in $\Trace{\cA}.$
\end{proof}

\noindent{\it Proof of Theorem \ref{p_main}.}
By Proposition \ref{idempotent} any $R_G$-irreducible invariant subspace $W$ of $\Q^n$ corresponds to
a minimal idempotent $E\in V_\Q(G)$ such that ${\sf Im}\{E\}=W$.
Furthermore,
since $\psi$ is an isomorphism between $V_\Q(G)$ and $S_\Q(G)$, the element
$\psi(E)$ is a minimal idempotent of $S_\Q(G)$ and therefore,
by Proposition \ref{p_idempotents}, $\psi(E) = \eps_d$ for some $d\in D(S_\Q(G))$.
Thus $W$ is $R_G$-irreducible invariant subspace of $\Q^n$ if and only if there exist
$d\in D(S_\Q(G))$
such that
\be \l{plm}
W = {\sf Im} \{\psi^{-1}(\eps_d)\}=
{\sf Im}\left\{\psi^{-1}(\sig_d)\Pi_{i=1}^\ell(I-\psi^{-1}(\sig_{f_i}))\right\}.
\ee

Observe now that if two idempotent matrices $A,$ $B$ commute
then for the matrix $C=AB=BA$
the equality
$${\sf Im}\{C\}= {\sf Im}\{A\}\cap{\sf Im}\{B\}$$ holds. Indeed, it is clear
that $${\sf Im}\{C\}\subseteq {\sf Im}\{A\}\cap{\sf Im}\{B\}.$$ On the other hand, if $z\in{\sf Im}\{A\}\cap{\sf Im}\{B\}$ then
$z=Ax=By$ for some vectors $x,y$ and \be \l{idem}Az=A(Ax)=Ax=z, \ \ \ Bz=B(By)=By=z.\ee It follows that $Cz=A(Bz)=Az=z$ and hence $z\in {\sf Im}\{C\}$.
Since Lemma \ref{corre} implies that $V_\Q(G)$ is commutative it follows now from \eqref{plm} that
$$
W={\sf Im}\left\{\psi^{-1}(\sig_d)\right\}\cap \left(\bigcap_{i=1}^\ell{\sf Im}\left\{(I-\psi^{-1}(\sig_{f_i}))\right\}\right).
$$

It was observed in the proof of Lemma \ref{p_D} that
${\sf Im}(\psi^{-1}(\sig_d))=V_{n/d}$.
Furthermore, since the image of any idempotent matrix consists of
its invariant vectors we have
${\sf Im}\{I -\psi^{-1}(\sig_d)\}={\sf Ker}\{\psi^{-1}(\sig_d)\}$. On the other hand,
since the matrix $\psi^{-1}(\sig_d)$ is symmetric,
${\sf Ker}\{\psi^{-1}(\sig_d)\}={\sf Im}\{\psi^{-1}(\sig_d)\}^\perp.$
Therefore,
$$
W = V_{n/d}\cap V_{n/f_1}^\perp\cap ... \cap V_{n/f_\ell}^\perp.
$$
Finally, Lemma \ref{p_D} implies that $n/d\in D(G)$ and
$n/f_1,...,n/f_\ell$ is a complete set of elements of $D(G)$ covered by $n/d$.
Hence, $W= U_{n/d}.$

\vskip 0.2cm
\noindent{\bf Remark.} If $G$ does not contain a full cycle, then Theorem \ref{p_main} fails to be true.
A simple example is provided by the group $S_5$ acting on two element subsets of
$\{1,2,3,4,5\}$. One can verify that in this way we obtain a primitive permutation group $G$ on $10$ points which
yields a permutation matrix representation $\rho_G$ of dimension $10$. However, the collection of $\rho_G$-invariant irreducible subspaces of $\Q^{10}$ is distinct from the collection  $U_1,U_{10}$ since  
$U_{10}$ is a direct sum
of two irreducible $\rho_G$-invariant subspaces of dimensions $4$ and $5$.

Notice also that Theorem \ref{p_main}
is not true for representations over $\C$.
In order to see this it is enough to take as $G$ any cyclic group.

\section{Description of $Q(z)$ satisfying $\phi_s(t)=0$}

\subsection{Geometry of $M_{P,a,b}$}
In notation of Section \ref{cycy} set $$W =V_{f_1}^\perp\cap ...\cap V_{f_\ell}^\perp,$$ where
$f_1,...,f_\ell$ is the set of all elements of $D(G_P)$
distinct from $n$.
Notice that since $n\in D(G_P)$ covers any other element of $D(G_P)$,
the subspace $W$ coincides with the subspace $U_n$ from
Theorem \ref{p_main} and therefore is $G_P$-invariant irreducible subspace of $\mathbb Q^n$.

Theorem \ref{p_main} together with Proposition \ref{mono} imply the following
important geometric property of $M_{P,a,b}.$

\bp \l{pp} The subspace $M_{P,a,b}$ contains the subspace $W.$
\ep

\pr Indeed, since by construction $M_{P,a,b}$ is a $G_P$-invariant subspace of $\Q^n$, it follows from
Theorem \ref{p_main} that either $M_{P,a,b}$ contains
$W$ or is orthogonal to $W.$ In the last case $M_{P,a,b}$ also would be orthogonal
to the complexification $W^{\C}$ of $W.$
Therefore, in order to prove the proposition
it is enough to find vectors
$\vec w\in W^{\C}$ and $\vec v\in M_{P,a,b}$
such that $(\vec v,\vec w)\neq 0.$

In order to find such $\vec w$ observe that the vectors
$$\vec w_i=(1,\varepsilon_n^j, \varepsilon_n^ {2j}, \,...\,,
\varepsilon_n^{(n-1)j}),$$ $1\leq  j \leq n,$
where $\varepsilon_n=exp(2\pi \sqrt{-1}/n),$
form an orthogonal basis of $\C^n.$ Furthermore, for $d\vert n$
vectors $\vec w_j$ for which $(n/d)\,\vert \,j$
form a basis of $V_d^{\C}.$ Therefore, the vector
$\vec w_1$ is orthogonal to $V_{f}^{\C}$
for any $f\in D(G_P),$ $f\neq n,$ and hence
$\vec w_1\in W^{\C}.$ Set $\vec w={\vec w}_1.$

Consider now two cases. Suppose first that $P(a)= P(b)$ and show that in this case
for the vector $\vec v\in M_{P,a,b}$ corresponding to equation \eqref{e1}
the inequality $(\vec v,\vec w)\neq 0$ holds.
Indeed, the equality $(\vec v,\vec w)= 0$ is equivalent to the equality
$$\sum_{s=1}^{d_a}\varepsilon_n^{a_s}/d_a =
\sum_{s=1}^{d_b}\varepsilon_n^{b_s}/d_b$$
which in its turn is equivalent to the statement
that the ``centers of mass'' of the sets $V(a)$ and $V(b)$
coincide. But this contradicts to Proposition \ref{mono} since
the center of mass of a system of points in $\C$
is inside of the convex envelope of this system and therefore
the centers of mass of disjointed sets must be distinct.

Similarly, if $P(a)\neq P(b)$ then $ (\vec v,\vec w)\neq 0$ for at least
one of two vectors corresponding to equations \eqref{e2}.
Indeed, otherwise
$$\sum_{s=1}^{d_a}\varepsilon_n^{a_s}/d_a =0, \ \ \ \ \
\sum_{s=1}^{d_b}\varepsilon_n^{b_s}/d_b=0$$ that
contradicts again to Proposition \ref{mono} since the fact that
the sets $V(a)$ and $V(b)$ are almost disjointed implies that at least
one of these sets is contained in an open half plane bounded by a line
passing through the origin and therefore has the center of mass
distinct from zero.
$\ \ \Box$

\subsection{Puiseux expansions of $Q(P^{-1}(z))$}
Let $\hat U\subset \C$ be a domain as in the proof of Proposition
\ref{mono}. Then, taking into account
our convention about the numeration of branches of
$P^{-1}(z)$, at points of $\hat U$ close enough to infinity
the function $Q(P^{-1}_i(z))$, $1\leq i \leq n,$
is represented by the converging series
\be \l{ps2}
Q(P^{-1}_i(z))=\sum_{k=-m}^{\infty}
s_k\varepsilon_n^{(i-1)k}z^{-\frac{k}{n}},
\ee
where $z^{\frac{1}{n}}$ denotes some fixed branch of the algebraic function
inverse to $z^n$ in $\hat U.$
Therefore, any relation of the form
\be \l{x}
\sum_{i=1}^{n}f_iQ(P^{-1}_i(z))=0, \ \ \ \ \ \ f_i\in \C,
\ee
is equivalent to the system
\be \l{vot}
\sum_{i=1}^{n}f_is_k\varepsilon^{k(i-1)}_n=0, \ \ \ k\geq -m.
\ee In particular, in view of Theorem \ref{t1}, the equality $\hat H(t)\equiv 0$ implies that 
for any $k\geq -m$ such that the coefficient $s_k$ of series \eqref{ps2} distinct from zero
the vector $\vec w_k$ is orthogonal to $M_{P,a,b}$.
This fact together with Proposition \ref{pp} imply the following
statement (cf. \cite{pp}, Theorem 4.1).

\bp \l{la} Let $Q(z)$ be a polynomial such that
$\hat H(t)\equiv 0.$ Then for any $k\geq -m$ such that the coefficient $s_k$ of series \eqref{ps2} is distinct from zero
there exists $f\in D(G_P),$ $f\neq n,$ such that
$(n/f)\,\vert\, k.$
\ep

\pr Indeed, if $s_k\neq 0$ then it follows from
\eqref{vot} that the vector $\vec w_k$ is orthogonal to $M_{P,a,b}^{\C}$ and
therefore by Proposition \ref{pp} is orthogonal to $W^{\C}$.
Since the subspace $(W^{\C})^{\perp}$ is generated by the vectors $\vec w_j,$ $(n/f)\,\vert \,j,$ $f\in D(G_P),$ $f\neq n,$ this implies that $\vec w_k$ is a linear combination of these vectors
and hence coincides with one of them since
the vectors $\vec w_i,$ $1\leq i \leq n,$ are linearly independent. Therefore, 
$(n/f)\,\vert \,k$ for some $f\in D(G_P),$ $f\neq n.$ \qed

For $f\in D(G_P)$, $f\neq n,$
set
$$\psi_{f}(z)=\sum_{\substack {k\geq -m \\ k\equiv 0\, \mod n/f}}s_{k} z^{-\frac{k}{n}},$$
where $s_k,$ $k\geq -m,$ are coefficients of series \eqref{ps2}. Clearly, $\psi_{f}(z)$ is
an analytic function in $\hat U.$

\bl \l{lape} For any $f\in D(G_P),$ $f\neq n,$ there exists $S_f(z)\in\C[z]$ such that \be \l{polk} \psi_{f}(z)=S_f(P_1^{-1}(z)).\ee
Furthermore, we have:
\be \l{xru} P(z)=A_1(B_1(z)), \ \ \ S_f(z)=R_1(B_1(z))\ee for some $A_1(z),B_1(z),R_1(z)\in \C[z]$ with $\deg B_1(z)>1.$
\el

\pr First, observe that since
$$1+(\varepsilon_n^{k})^f+(\varepsilon_n^{k})^{2f}+\dots +(\varepsilon_n^k)^{n-f}$$
equals $n/f$ if $n\vert(fk)$ and zero otherwise, it follows from \eqref{ps2} that
the equality
\be \l{ra}\left(\frac{n}{f}\right)\psi_{f}(z)=
Q(P^{-1}_{1}(z))+
Q(P^{-1}_{f+1}(z))+Q(P^{-1}_{2f+1}(z))
+ ... + Q(P^{-1}_{n-f+1}(z))
\ee holds.

Let now $\Omega_P$ be a field generated by all branches of $P^{-1}(z)$ 
considered as elements of some fixed algebraic closure of $\C(z)$. Recall that
the Galois group of the extension $[\Omega_P:\C(z)]$
is permutation equivalent to the group $G_P$ and  
under the Galois correspondence
to the stabilizer of $P^{-1}_1(z)$ in $G_P$ corresponds the invariant subfield $\C(P^{-1}_1(z))$ of $\Omega_P$.
Since $f\in D(G_P),$
the collection of branches appearing in the right part of equality \eqref{ra} is a block of an imprimitivity system of
$G_P$ containing $P^{-1}_1(z)$. Therefore, equality \eqref{ra} implies that
the function $\psi_{f}(z)\in \Omega_P$ is invariant with respect
to the action of the stabilizer of $P^{-1}_1(z)$ in $G_P$ and hence
is contained in the field $\C(P^{-1}_1(z))$.
So, there exists a rational function $S_f(z)$ such that equality \eqref{polk} holds.
Furthermore, since the analytic continuation of the right side of \eqref{ra}
has no poles in $\C$ the function $S_f(z)$ is a polynomial.
Finally, since branches appearing in the right part of equality \eqref{ra} form a block,
it is easy to see that 
$$S_f(P_1^{-1}(z))=S_f(P_{lf+1}^{-1}(z)), \ \ \ 1\leq l \leq n/f-1,$$ 
and hence the last part of the lemma follows from Lemma \ref{lll}.
\qed

\subsection{Proof of Theorem 1.1} In view of Theorem \ref{t1} we essentially must 
show that the conclusion of the theorem holds for any non zero polynomial $Q(z)$ such that 
$\hat H(t)\equiv 0.$ So, abusing the notation, below we will mean by a solution of the polynomial moment problem such a polynomial $Q(z)$.
The proof is by induction on the number $i(P)$ of
imprimitivity systems of the group $G_P$.
If $i(P)=2$, that is if $G_P$ has only trivial imprimitivity systems, then Proposition \ref{la} implies that for any non-zero coefficient $s_j,$ $j\geq m,$ of  \eqref{ps2} the number $k$ is a multiple of $n$.
Therefore, all the functions $Q(P^{-1}_i(z)),$ $1\leq i \leq n,$ are equal between themselves and hence $Q(z)=R(P(z))$ for some polynomial $R(z)$ by 
Lemma \ref{lll}. Furthermore, necessarily
$P(a)=P(b)$. Indeed,
otherwise after the change of variable $z=P(z)$ we would
obtain that the polynomial $R(z)$ is orthogonal to all powers
of $z$ on the segment  $[P(a),P(b)]$. However, for
$$P(z)=z, \ \ Q(z)=R(z), \ \ a=P(a), \ \ b=P(b)$$
any of relations \eqref{e2} reduces to the
equality $R(z)\equiv 0$ in contradiction with the condition $Q(z)\not\equiv 0$
(of coarse instead of Proposition 2.1 we also
could use the Weierstrass theorem).
Therefore,
if $i(P)=2$ then all solutions of the polynomial moment problem for $P(z)$ are reducible
(cf. \cite{pa2}, Theorem 1 and \cite{pp}, Theorem 5.3).

Suppose now that the theorem is proved for all $P(z)$ with
$i(P)<l$ and let $Q(z)$ be a non-zero solution of the polynomial moment problem for a polynomial $P(z)$
of degree $n$ with $i(P)=l.$
If $Q(z)=R(P(z))$ for some polynomial $R(z)$ then one can show as above that
$P(a)=P(b)$ and 
$Q(z)$ is reducible.
Otherwise there exists a non-zero coefficient $s_{j_1},$ $j_1\geq m,$ of expansion
\eqref{ps2} such that $j_1$ is not a multiple of $n$. By Proposition \ref{la} this implies that
there exists $f_1\in D(G_P),$ $f_1\neq n,$ such that
$(n/f_1) \,\vert\, j_1$.
Furthermore, by Lemma \ref{lape} there exists a polynomial
$S_1(z)$ such that
$\psi_{f_1}(z)=S_1(P^{-1}_1(z))$ and equalities
$$ P(z)=A_1(B_1(z)), \ \ \ S_1(z)=R_1(B_1(z))$$ hold for some $A_1(z),B_1(z),R_1(z)\in \C[z]$ with $\deg B_1(z)>1.$

Define a polynomial $T_1(z)$ by the equality $T_1(z)=Q(z)-S_1(z)$. Then for any $i,$ $1\leq i \leq n,$ we have:
$$Q(P^{-1}_i(z))=S_1(P^{-1}_i(z))+T_1(P^{-1}_i(z)).$$
Since by
construction the intersection of the supports of the series
$S_1(P^{-1}(z))$ and $T_1(P^{-1}(z))$ is empty, if the series $Q(P^{-1}_i(z)),$ $1\leq i \leq n,$ satisfy some linear relation over $\C$ then the series $S_1(P^{-1}_i(z)),$ $1\leq i \leq n,$ and $T_1(P^{-1}_i(z)),$
$1\leq i \leq n,$ also satisfy this relation. It follows now from Theorem
\ref{t1} that
each of germs defined in a neighborhood of infinity by the integrals
$$
\hat H_1(t)= \int_{\Gamma_{a,b}} \frac{S_1(z)P^{\prime}(z)dz}{P(z)-t}\,, \ \ \ \
\hat F_1(t)= \int_{\Gamma_{a,b}} \frac{ T_1(z)P^{\prime}(z)dz}{P(z)-t}\,,
$$
vanishes or in other words the polynomials $S_1(z)$ and $R_1(z)$ are solutions of the polynomial moment problem for $P(z).$
Moreover, by
construction the Puiseux series of $T_1(P^{-1}(z))$
contains no non-zero coefficients with indices which are multiple of
$n/ f_1$. In particular, this implies that all coefficients of $T_1(P^{-1}(z))$ whose indices are multiple of $n$
vanish and hence
$T_1(z)$ may not have the form $T_1(z)=R(P(z))$ for some $R(z)\in \C[z]$ unless $T_1(z)\equiv 0.$

If $T_1(t)\neq 0$ then arguing as above we conclude that
there exist $f_2\in D(G_P),$ $f_2\neq f_1,$ $f_2\neq n$, and
polynomials $S_2(z), T_2(z), R_2(z), A_2(z), B_2(z) \in \C[z]$ with
$\deg B_2(z)>1$ such that the following conditions hold:
$$T_1(P^{-1}(z))=S_2(P^{-1}(z))+T_2(P^{-1}(z)),$$
$$ P(z)=A_2(B_2(z)), \ \ \ S_2(z)=R_2(B_2(z)),$$
the germs
$$
\hat H_2(t)= \int_{\Gamma_{a,b}} \frac{S_2(z)P^{\prime}(z)dz}{P(z)-t}\,, \ \ \ \
\hat F_2(t)= \int_{\Gamma_{a,b}} \frac{T_2(z)P^{\prime}(z)dz}{P(z)-t}\,
$$ vanish,
and the Puiseux expansion of $T_2(P^{-1}(z))$
contains no non-zero coefficients whose indices are multiple of
$n/ f_1$ or $n/ f_2$.

Since the number of divisors of $n$ is finite, continuing in this way, after a finite number of steps we will arrive to
a decomposition of the function $Q(z)$ into a sum of polynomials $S_s(z),$ $1\leq s \leq r,$
$$Q(z)=S_1(z)+S_2(z)+\dots +S_r(z)$$ such that
the germs
$$
\hat H_s(t)= \int_{\Gamma_{a,b}} \frac{S_s(z)P^{\prime}(z)dz}{P(z)-t}\,, \ \ \ \
1 \leq s \leq r,
$$
vanish and $$P(z)=A_s(B_s(z)), \ \ \ \
S_s(z)=R_s(B_s(z)), \ \ \ \ 1 \leq s \leq r, $$
for some $R_s(z), A_s(z), B_s(z) \in \C[z]$ with $\deg B_s(z)>1.$

In order to finish the proof it is enough to show any polynomial $S(z)$ from the collection $S_s(z),$ $1\leq s \leq r,$ is a sum of reducible solutions
of the polynomial moment problem for $P(z)$. So, take some $S(z)$ and let
$R(z), A(z), B(z),$ $\deg B(z)> 1,$ be polynomials
such that
$$P(z)=A(B(z)), \ \ \ \
S(z)=R(B(z)). $$
If $B(a)=B(b)$ then $S(z)$ itself is a reducible
solution. Otherwise, since
$$\int_{\Gamma_{a,b}} \frac{S(z)P^{\prime}(z)dz}{P(z)-t}
= \int_{B(\Gamma_{a,b})}\frac{R(z)A^{\prime}(z)dz}{A(z)-t},$$ we conclude that
the polynomial
$R(z)$ is a solution
of the polynomial moment problem for the polynomial $A(z)$ (and the points $B(a),$ $B(b)$). Since the condition
$\deg B(z)>1$ implies that $i(A)<i(P)$
it follows from the induction assumption that
there exist polynomials $V_{1}(z),$ $V_{2}(z), \dots , V_{j}(z)$ such that
$$R(z)= V_{1}(z)+V_{2}(z)+ \dots + V_{j}(z)$$
and $$V_{e}(z)=\tilde V_{e}(U_{e}(z)),\ \ \
A(z)=\tilde A_{e}(U_{e}(z)), \ \ \  U_{e}(B(a))=U_{e}(B(b)),$$
for some  $\tilde V_{e}(z),\tilde A_{e}(z), U_{e}(z) \in \C[z],$ $1 \leq e\leq j.$

Set now $$E_e(x)=V_e(B(x)), \ \ \ W_e(z)=U_{e}(B(z)),\ \ \ 1 \leq e\leq j.$$ Then 
$$S(z)=E_1(z)+E_2(z)+\dots +E_j(z),$$ where 
for each
$e,$ $1\leq e \leq j,$ we have:
$$E_e(z)= \tilde V_{e}(W_e(z)),\ \ \ P(z)=\tilde A_e(W_e(z)), \ \ \ W_e(a)=W_e(b).$$
Therefore,  
$S(z)$ is a sum of reducible solutions. \qed

\vskip 0.2cm
\noindent{\bf Remark.} Theorem 1.1 implies that if for a given polynomial $P(z)$ the corresponding polynomial moment problem has non-reducible solutions, then $P(z)$ has at least one ``double decomposition''
$$P=A\circ B=C\circ D$$
such that $$B(z)\notin \C(D(z)), \ \ \ D(z)\notin \C(B(z)).$$ 
Notice that this condition
is quite restrictive. Namely, the results of Engstrom \cite{en} and Ritt \cite{ri}
imply that if polynomials $A,B,$ $C,D$ satisfy the equation
$$ A\circ B= C\circ D $$
then there exist polynomials
$\hat A, \hat B, \hat C, \hat D, U, V$ such that
$$A=U\circ \hat A, \ \  C=U\circ \hat C,\ \ B=\hat B\circ V, \ \  D=\hat D \circ V, \ \ \hat A\circ \hat B=\hat C\circ \hat D,$$ and
up to a possible replacement of
$\hat A$ by $\hat C$ and $\hat B$ by $\hat D$ either
$$\hat A\circ \hat B\sim z^n \circ z^rR(z^n),  \ \ \ \ \ \ \hat C\circ \hat D
\sim  z^rR^n(z) \circ z^n,$$
where $R(z)$ is a polynomial, $r\geq 0,$ $n\geq 1,$ and
$\GCD(n,r)=1,$ or $$\hat A\circ \hat B \sim T_n \circ T_m, \ \ \ \ \ \ \hat C\circ \hat D\sim T_m \circ T_n,$$
where $T_n(z),T_m(z)$ are the corresponding Chebyshev polynomials, $n,m\geq 1,$ and $\GCD(n,m)=1.$

Notice however that a polynomial $P(z)$ may have more than {\it one} 
double decomposition satisfying the condition above.
Indeed, for example for any distinct prime divisors $p_1,p_2$ of a number $n$
we have
$$T_n(z)=T_{n/p_1}(T_{p_1}(z))=T_{n/p_2}(T_{p_2}(z))$$ and 
$$T_{p_1}(z)\notin \C(T_{p_2}(z)), \ \ \ T_{p_2}(z)\notin \C(T_{p_1}(z).$$
It would be interesting to investigate what
conditions should be imposed on the collection $P(z),a,b$
in order to conclude that any solution of the polynomial moment problem for $P(z)$ can be represented as a sum
{\it at most} $r$ reducible solutions, where $r\geq 1$ is a fixed number.

\bibliographystyle{amsplain}

\end{document}